\documentclass[11pt]{article}
\usepackage{amssymb}
\usepackage{mathrsfs}
\usepackage{latexsym,amsmath,amssymb,amsfonts,epsfig,graphicx,cite,psfrag}
\usepackage{eepic,color,colordvi,amscd}
\usepackage{comment}

\newtheorem{theo}{Theorem}[section]

\newtheorem{lem}[theo]{Lemma}
\newtheorem{coro}[theo]{Corollary}

\newtheorem{prop}[theo]{Proposition}

\def\qed{\hfill \rule{4pt}{7pt}}

\def\pf{\noindent {\it Proof. }}

\begin{document}

\title{The Kelmans-Seymour conjecture I:  special separations}

\author{Dawei He\footnote{dhe9@math.gatech.edu; Partially supported by NSF grant through X. Yu},
Yan Wang\footnote{yanwang@gatech.edu; Partially supported by NSF grant through X. Yu},  
Xingxing Yu\footnote{yu@math.gatech.edu; Partially supported by NSF
  grants DMS-1265564 and CNS-1443894} \\
\bigskip \\
School of Mathematics\\
Georgia Institute of Technology\\
Atlanta, GA 30332}


\date{}

\maketitle

\begin{abstract}
 Seymour and, independently, Kelmans 
conjectured in the 1970s that every 5-connected nonplanar graph contains a
subdivision of $K_5$. This conjecture was proved by Ma and Yu for graphs containing $K_4^-$, and  
an important step in their proof is to deal with a 5-separation in the graph with a planar side. 
In order to establish the Kelmans-Seymour conjecture for all graphs, 
we need to consider 5-separations and 6-separations with less restrictive structures.
The goal of this paper is to deal with special 5-separations and 6-separations, including those with an apex side. 
Results will be used in subsequent papers to prove the Kelmans-Seymour conjecture.


\end{abstract}

\newpage 

\section{Introduction}
Let $K$ be a graph; we use $TK$ to denote a subdivision of $K$ and call the vertices of the $TK$ corresponding to the vertices of  
$K$ its {\it branch} vertices. 
Kuratowski's theorem states that a graph is planar iff it contains neither $TK_{3,3}$ nor $TK_5$. Graphs containing no 
$TK_{3,3}$ can be constructed from planar graphs and copies of $K_5$ by pasting them along 
cliques of size at most two. The structure of graphs containing no $TK_5$ is not well understood. 
Kelmans~\cite{Ke79} and, independently, Seymour~\cite{Se77} conjectured that 5-connected nonplanar graphs 
must contain $TK_5$. Thus, if the Kelmans-Seymour conjecture is true then  graphs containing no $TK_5$ is planar or 
admits a cut of size at most 4. Note that the requirement on connectivity is best possible (see, for example, 
$K_{4,4}$). 

Ma and Yu \cite{MY10, MY13} proved the Kelmans-Seymour conjecture for graphs containing $K_4^-$, and  Kawarabayashi, Ma and Yu~\cite{KMY15} 
proved the Kelmans-Seymour conjecture for graphs containing $K_{2,3}$. We refer the reader 
to \cite{MY10,MY13,KMY15} for problems and results (as well as references) related to the Kelmans-Seymour conjecture. 

 It turns out that $K_4^-$ is the right intermediate structure for studying the Kelmans-Seymour structure. 
By a result of Kawarabayashi \cite{Ka02}, any 5-connected graph containing   no $K_4^-$ has an edge $e$ that is {\it contractible} 
(i.e.,  $G/e$ is also 5-connected). Thus, our strategy for proving the Kelmans-Seymour conjecture is to keep contracting edges 
incident with a special vertex to produce a smaller 5-connected graph. To avoid trivial components associated with $5$-cuts or $6$-cuts, 
we also contract triangles (but we give preference to edges). 

We now give a more detailed description of our strategy. For a graph 
$G$ and a connected subgraph $M$ of $G$ (respectively, an edge $e$ of $G$), we use $G/M$ (respectively, $G/e$) to denote the graph obtained from $G$ 
by contracting $M$ (respectively, $e$). 
Let $G$ be a 5-connected nonplanar graph containing no $K_4^-$. Then $G$ contains an edge $e$ such that $G/e$ is 5-connected. If $G/e$ is planar, 
we can apply a discharging argument. So assume $G/e$ is not planar. Let $M$ be a maximal connected subgraph of $G$ such that $G/M$ is 5-connected and nonplanar. 
Let $z$ denote the vertex representing the contraction of $M$, and let $H=G/M$. Then one of the following holds.
\begin{itemize}
\item [(a)] $H$ contains a $K_4^-$ in which  $z$ is of degree 2.
\item [(b)] $H$ contains a $K_4^-$ in which  $z$ is of degree 3.
\item [(c)] $H$ does not contain $K_4^-$, and there exists $T\subseteq H$ such that $z\in V(T)$, $T\cong K_2$ or $T\cong K_3$, and $H/T$ is 5-connected and planar.
\item [(d)] $H$ does not contain  $K_4^-$, and for any $T\subseteq H$ with $z\in V(T)$ and $T\cong K_2$ or $T\cong K_3$, $H/T$ is not 5-connected.
\end{itemize}
We plan a series of five papers to establish the Kelmans-Seymour conjecture.  
The purpose of this paper is  to prove several results about 5-separations and 6-separations, which will be used in subsequent papers to 
reduce (c) and (d) to (a) or (b). Note that 5-separations and 6-separations arise naturally when (d) occurs. 
In the second paper, we will handle (a). 
The third and fourth papers take care of (b). 
In the final paper, we will deal with (c) and (d). In our arguments throughout this series, we frequently encounter the  
case when $H-z$ contains $K_4^-$; so the exclusion of $K_4^-$ (in \cite{MY10, MY13}) is very useful.

One of the main steps in the proofs in \cite{MY10,MY13} is to deal with 5-connected nonplanar graphs that admit a
5-separation with a planar side. A {\it separation} in a graph $G$ consists of a pair of subgraphs $G_1, G_2$ of $G$, denoted as $(G_1,G_2)$, such that
$E(G_1) \cup E(G_2)=E(G)$, $E(G_1\cap G_2)=\emptyset$, and neither $G_1$ nor $G_2$ is a subgraph of the other. 
The {\it order} of this separation is $|V(G_1\cap G_2)|$, and $(G_1,G_2)$ is said to be a {\it $k$-separation} if
its order is $k$.  Let $G$ be a graph and $A\subseteq V(G)$. We say that $(G,A)$ is {\it planar} if $G$ has a plane representation in which 
the vertices in $A$ are incident with a common face.
Ma and Yu \cite{MY10} proved that if $G$ has a 5-separation $(G_1,G_2)$ such that $|V(G_i)|\ge 7$ (for $i=1,2$) and 
$(G_2,V(G_1\cap G_2))$ is planar then $G$ contains a $TK_5$ with all branch vertices contained in $G_2$. 
In order  to establish the Kelmans-Seymour conjecture for graphs containing no $K_4^-$, 
we need to study 5-separations and 6-separations with less restrictive structures.

Let $G$ be a graph and $a\in V(G)$. Then $G-a$ denotes the subgraph obtained from $G$ by deleting $a$ and all edges of $G$ incident with $a$. 
For a positive integer $n$, we let $[n]=\{1,\ldots, n\}$. We often represent a path (or cycle) by a sequence of vertices.  
The following result deals with one type of 5-separations.

\begin{theo}\label{apexside1}
Let $G$ be \textsl{a} $5$-connected nonplanar graph and let $(G_1, G_2)$ be \textsl{a} $5$-separation in $G$. 
Suppose $|V(G_i)|\geq 7$ for $i\in [2]$, $a\in V(G_1\cap G_2)$, and $(G_2-a,V(G_1\cap G_2)-\{a\})$ is planar. Then one of the following holds.
\begin{itemize}
\item [$(i)$]  $G$ contains \textsl{a} $TK_5$ in which $a$ is not \textsl{a} branch vertex.
\item [$(ii)$] $G-a$ contains $K_4^-$.
\item [$(iii)$] $G$ has \textsl{a} $5$-separation $(G_1',G_2')$ such that $V(G_1'\cap G_2')=\{a, a_1,a_2,a_3,a_4\}$ and $G_2'$ is 
the graph obtained from the edge-disjoint union of the $8$-cycle $a_1b_1a_2b_2a_3b_3a_4b_4\allowbreak a_1$ 
and the  $4$-cycle $b_1b_2b_3b_4b_1$ by adding $a$ and the edges $ab_i$ for $i\in [4]$.
\end{itemize}
\end{theo}

Let $G$ be a graph. For $S\subseteq V(G)$, we use $G[S]$ to denote the subgraph of $G$ induced by $S$, and let $G-S=G[V(G)-S]$. 
For $H\subseteq G$, we write $G[H]$ instead of $G[V(H)]$.  For $S\subseteq E(G)$, $G-S$ denotes the graph obtained from $G$ by 
deleting the edges in $S$.  Another type of 5-separations considered in this paper are those  $(G_1,G_2)$
with the property that $G[V(G_1\cap G_2)]$ contains a triangle.

\begin{theo}
\label{5cut_triangle}
Let $G$ be \textsl{a} $5$-connected graph and $(G_1,G_2)$ be \textsl{a} $5$-separation in $G$. Suppose that $|V(G_i)|\ge 7$ for $i\in [2]$ and 
$G[V(G_1\cap G_2)]$ contains \textsl{a} triangle $aa_1a_2a$. Then one of the following holds. 
\begin{itemize}
\item [$(i)$]  $G$ contains \textsl{a} $TK_5$ in which $a$ is not \textsl{a} branch vertex.
\item [$(ii)$] $G-a$ contains $K_4^-$.
\item [$(iii)$] $G$ has \textsl{a} $5$-separation $(G_1',G_2')$ such that $V(G_1'\cap G_2')=\{a, a_1,a_2,a_3,a_4\}$ and $G_2'$ is 
the graph obtained from the edge-disjoint union of the $8$-cycle $a_1b_1a_2b_2a_3b_3a_4b_4\allowbreak a_1$ 
and the $4$-cycle $b_1b_2b_3b_4b_1$ by adding $a$ and the edges $ab_i$ for $i\in [4]$.
\item [$(iv)$] For any distinct $u_1,u_2,u_3\in N(a)-\{a_1,a_2\}$, $G-\{av: v \not\in \{a_1,a_2,u_1,u_2,u_3\}\}$ contains $TK_5$. 
\end{itemize}
\end{theo}

 In the applications of Theorems~\ref{apexside1} and \ref{5cut_triangle}, the vertex $a$ will represent the special vertex 
resulting from the contraction of a connected subgraph formed by a sequence of edges and triangles. (Allowing the contraction of triangles will ensure that $|V(G_i)|\ge 7$ 
in the applications of Theorems~\ref{apexside1} and \ref{5cut_triangle}.) 
So $(i)$ of Theorems~\ref{apexside1} and \ref{5cut_triangle} gives a $TK_5$ in the original graph, and $(ii)$ of 
Theorems~\ref{apexside1} and \ref{5cut_triangle} allows us to apply the result of Ma and Yu~\cite{MY13} to get a $TK_5$ in the original graph.  The $TK_5$ given in  $(iv)$ of 
Theorem~\ref{5cut_triangle} can be used to derive a $TK_5$ in the original graph. 
 When $(iii)$ of Theorems~\ref{apexside1} and \ref{5cut_triangle} occurs, we will use Proposition~\ref{8cycle} below, whose proof is included in this section (as it is short).

Let $G$ be a graph. By $H\subseteq G$ we mean that $H$ is a subgraph of $G$. When $K\subseteq G$ and $L\subseteq G$, we let $K-L=K-V(K\cap L)$.  
For $S\subseteq V(G)$, we may view  $S$ as a subgraph of $G$ with vertex set $S$ and edge set $\emptyset$. 
For $H\subseteq G$, $N_G(H)$ denotes the neighborhood of $H$ (not 
including the vertices in $V(H)$).  For any $x\in V(G)$, we use $N_G(x)$ to denote the neighborhood of $x$ in $G$. 
When understood, the reference to $G$ may be dropped. 
We may view paths as sequences of vertices. The {\it ends} of
a path $P$ are the vertices of the minimum degree in $P$, and all
other vertices of $P$ (if any) are its {\it
  internal} vertices. A collection of paths are said to be {\it independent} if no vertex
of any path in this collection is an internal vertex of any other path in the collection.

\begin{prop}
\label{8cycle}
Let $G$ be \textsl{a} 5-connected nonplanar graph, $(G_1,G_2)$ \textsl{a} 5-separation in $G$, $V(G_1 \cap G_2)=\{a,a_1,a_2,a_3,a_4\}$ such that  
$G_2$ is the graph obtained from the edge-disjoint union of the $8$-cycle $a_1b_1a_2b_2a_3b_3a_4b_4a_1$ 
and the $4$-cycle $b_1b_2b_3b_4b_1$ by adding $a$ and the edges $ab_i$ for $i\in [4]$. Suppose $|V(G_1)|\ge 7$. 
Then, for any $u_1, u_2 \in N(a) - \{b_1,b_2,b_3\}$, $G - \{av: v \not\in \{b_1,b_2,b_3,u_1,u_2\} \}$ contains  $TK_5$.
\end{prop}

\pf
By symmetry between $u_1$ and $u_2$, we may assume that $u_1 \neq b_4$. For convenience, let 
$G':=G - \{av: v \not\in \{b_1,b_2,b_3,u_1,u_2\} \}$.

First, suppose $u_1 \in V(G_1 - G_2)$. 
Since $G$ is $5$-connected, $G_1-a$ has independent paths $Q_1,Q_2,Q_3,Q_4$ from $u_1$ to $a_1,a_2,a_3,a_4$, respectively. 
Then $G[\{a,b_1,b_2,b_3,u_1\}] \cup (Q_1 \cup a_1b_1) \cup (Q_2 \cup a_2b_2) \cup (Q_3 \cup a_3b_3) \cup b_1b_4b_3$ is a $TK_5$ in $G'$
with branch vertices $a, b_1,b_2,b_3,u_1$.

Now suppose $u_1\in \{a_2,a_3\}$. By symmetry, we may assume $u_1 = a_2$. 
Since $G$ is $5$-connected, $G_1-a$ contains a path $R$ from $a_2$ to $a_3$. 
Then $G[\{a,b_1,b_2,b_3,u_1\}] \cup (R \cup a_3b_3) \cup b_1b_4b_3$ is a $TK_5$ in $G'$ with branch vertices $a, b_1,b_2,b_3,u_1$. 

Finally, assume $u_1\in \{a_1,a_4\}$. By symmetry, we may assume  $u_1 = a_1$. 
If $G_1 -a$ has independent paths $R_1,R_2$ from $a_1$ to $a_2,a_3$, respectively, 
then $G[\{a,b_1,b_2,b_3,u_1\}] \cup (R_1 \cup a_2b_2) \cup (R_2 \cup a_3b_3) \cup b_1b_4b_3$ is a $TK_5$ in $G'$ with branch vertices $a, b_1,b_2,b_3,u_1$. 
So we may assume that such $R_1,R_2$ do not exist. Then there exists a separation $(K,L)$ in $G_1 - a$ such that 
$|V(K\cap L)|\le 1$, $a_1 \in V(K)$, and $\{a_2,a_3\} \subseteq V(L)$. 
Since $G$ is $5$-connected, $V(K)=\{a_1,a_4\}\cup V(K\cap L)$ and $V(L)=\{a_2,a_3\}\cup V(K\cap L)$. But this implies 
$|V(G_1)|\le 6$, a contradiction. \qed

\medskip

The next result  deals with a special type of  $6$-separations, whose proof makes use of Theorems~\ref{apexside1} and \ref{5cut_triangle}.
\begin{theo} \label{N(x1)capA}
Let $G$ be \textsl{a} $5$-connected graph and $x\in V(G)$, and let 
$(G_1,G_2)$ be \textsl{a} $6$-separation in $G$ such that 
$x\in V(G_1\cap G_2)$,  $G[V(G_1\cap G_2)]$ contains a triangle
$xx_1x_2x$, $|V(G_i)|\ge 7$ for $i\in [2]$. Moreover, assume that
$(G_1,G_2)$ is chosen so that, subject to $\{x,x_1,x_2\} \subseteq
V(G_1\cap G_2)$ and $|V(G_i)|\ge 7$ for $i\in [2]$, $G_1$ is minimal.  
Let $V(G_1 \cap G_2) = \{x,x_1,x_2,v_1,v_2,v_3\}$. 
Then $N(x)\cap V(G_1-G_2)\ne \emptyset$, or  one of the following holds. 
\begin{itemize}
\item[$(i)$] $G$ contains \textsl{a} $TK_5$ in which $x$ is not \textsl{a} branch vertex.
\item[$(ii)$] $G$ contains $K_4^-$. 
\item[$(iii)$] There exists $x_3\in N(x)$ such that for any distinct $y_1,y_2\in N(x)-\{x_1,x_2,x_3\}$, 
$G-\{xv:v\notin \{x_1,x_2,x_3,y_1,y_2\}\}$ contains $TK_5$.
\item [$(iv)$] For some $i\in [2]$ and some $j \in [3]$, 
$N(x_i)\subseteq V(G_1-G_2)\cup \{x,x_{3-i}\}$, and  any three
independent paths in $G_1-x$ from $\{x_1,x_2\}$ to $v_1,v_2,v_3$, respectively, 
with two from $x_i$ and one from $x_{3-i}$, 
must contain \textsl{a} path from $x_{3-i}$ to $v_j$.
\end{itemize}
\end{theo}

Note that in $(ii)$ of Theorem~\ref{N(x1)capA}, we ask that $G$ contain $K_4^-$
instead of  the stronger statement ``$G-a$ contains $K_4^-$'' as in (ii) of
Theorems~\ref{apexside1} and \ref{5cut_triangle}. Thus, $(iii)$ of
Theorems~\ref{apexside1} and \ref{5cut_triangle} does not occur in
Theorem~\ref{N(x1)capA} as it gives $(ii)$ of Theorem~\ref{N(x1)capA}.

\medskip

This paper is organized as follows. In Section 2, we use the discharging technique to prove two lemmas about $K_4^-$ in apex graphs. 
(A graph is apex if it has a vertex whose removal results in a planar graph.)
In Section 3, we collect a number of known results and prove Theorem~\ref{apexside1}.
In Section 4, we prove a result about apex graphs (from which one can see how (c) might be taken care of).
In Section 5, we prove Theorem~\ref{5cut_triangle}. In Section 6, we prove Theorem~\ref{N(x1)capA}, using 
Theorems~\ref{apexside1} and \ref{5cut_triangle}.

We end this section with additional notation and terminology. Let $G$ be a graph. 
Let $K\subseteq G$, $S\subseteq V(G)$, and $T$ a collection of 2-element subsets of $V(K)\cup S$; 
then $K+(S\cup T)$ denotes the graph with vertex set $V(K)\cup S$ and edge set $E(K)\cup T$. If
$T=\{\{x,y\}\}$ and $x,y\in V(K)$,  we write $K+xy$ instead of $K+\{\{x,y\}\}$.

 A  set $S\subseteq V(G)$ is a {\it $k$-cut} (or a {\it cut} of size $k$)
in a graph $G$, where $k$ is a positive integer,  if $|S|=k$ and $G$ has a
separation $(G_1,G_2)$ such that $V(G_1\cap G_2)=S$ and $V(G_i-S)\ne
\emptyset$ for $i\in [2]$. If $v\in V(G)$ and $\{v\}$ is a cut
of $G$, then $v$ is said to be a {\it cut vertex} of $G$.

Given a path $P$ in a graph and $x,y \in V(P)$, $xPy$ denotes the
subpath of $P$ between $x$ and $y$ (inclusive).
 A path $P$ with ends $u$ and $v$ (or an $u$-$v$ path) is also said to be {\it
from $u$ to $v$} or {\it between $u$ and $v$}.

A {\it plane} graph is a graph drawn in the plane with no edge-crossings. The unbounded face of a plane graph is usually called its {\it 
outer face}. The boundary of the outer face of a connected plane graph is the {\it outer walk} of the graph (or   {\it outer cycle} when it is a cycle). 
It is well known that if $G$ is a 2-connected plane graph then every facial boundary of $G$ is a cycle. 
Let  $D$ be a cycle in a plane graph. Given $x,y \in V(D)$, if $x\ne y$ then $xDy$ denotes the subpath of $D$ between $x$ and $y$ (inclusive) in clockwise order; 
and if $x=y$ then $xDy$ is simply the trivial path with the single vertex $x=y$.

\section{Discharging and $K_4^-$}

In this section, we prove results about $K_4^-$ in certain apex graphs,  using the discharging technique. First, we give a simple lemma
on discharging. For a plane graph $G$, let $F(G)$ denote the set of all faces of $G$ and, for each $f\in F(G)$, let $d(f)$ denote the 
number of edges of $G$ incident with $f$ (with each cut edge counted twice).  

\begin{lem}\label{discharge}
Let $G$ be a connected plane graph and let $\sigma: V(G)\cup F(G) \longrightarrow \mathbb{Z}^+$, the set of nonnegative integers,  
such that $\sigma(t)=4-d(t)$ for all $t\in V(G)\cup F(G)$. Let $\tau$ be obtained from $\sigma$ as follows: For each 
$f\in F(G)$ with $d(f)=3$, choose two vertices incident with $f$ and send charge $1/2$ from $f$ to each of these two vertices. Then 
$$\sum\limits_{v\in V(G)}\sigma(v) + \sum\limits_{f\in F(G)}\sigma(f) = 8,$$
and if $K_4^-\not\subseteq G$ then, for $v\in V(G)$, 
$$
\tau(v) \leq \left\{ \begin{array}{ll}
 4-3k, &\mbox{ if $d(v)=4k$;} \\
 3-3k, &\mbox{ if $d(v)=4k+1$;}\\
 5/2-3k, &\mbox{ if $d(v)=4k+2$;}\\
 3/2-3k, &\mbox{ if $d(v)=4k+3$.}      
 \end{array} \right.
$$
\end{lem}

\pf  By Euler's formula, we have
$$\sum\limits_{v\in V(G)}\sigma(v) + \sum\limits_{f\in F(G)}\sigma(f) = 8.$$

Now suppose $K_4^-\not\subseteq G$. Then for each $v\in V(G)$, $v$ is contained in at most $\lfloor d(v)/2\rfloor$ 
facial triangles. So $\tau(v)\leq \sigma(v) + \lfloor d(v)/2\rfloor/2 = 4-d(v) + \lfloor d(v)/2\rfloor/2$. By considering 
$d(v)$ modulo 4, we get the desired bounds on $\tau(v)$. 
\qed

\medskip

To state the remaining results in this section, we need a concept on connectivity.
Let $G$ be a graph and $A\subseteq V(G)$, and let $k$ be a positive integer.  We say that 
$G$ is {\it $(k,A)$-connected} for some positive integer $k$ if, for any cut $T$ of $G$ with $|T|< k$, each component of 
$G-T$ contains a vertex from $A$.  Recall that $(G,A)$ is planar if $G$ has a plane representation in which the vertices in $A$ are incident with 
a common face.

\begin{lem}\label{planarside}
Let $G$ be \textsl{a} connected graph and $A\subseteq V(G)$ such that $|A|=5$, $|V(G)|\ge 7$, $G$ is $(5,A)$-connected, and $(G,A)$ is planar. 
Then, for any $a\in A$, $G-a$ contains $K_4^-$. 
\end{lem} 

\pf It suffices to prove the lemma for the case when $A$ is an independent set in $G$. So we assume that $A$ is an independent set in $G$. 
First, we show that $G-a$ is connected for any $a\in A$. For, if not, then we may  
let $C_1, C_2$ be two components of $G-a$ such that $V(C_1)\not\subseteq A$. 
Since $G$ is $(5,A)$-connected, $V(C_i)\cap A\ne \emptyset$ for $i\in [2]$. 
 Now $S:=(A\cap V(C_1))\cup \{a\}$ is a cut in $G$ such that  $|S|\le 4$ and $G-S$ has a component contained in $C_1-A$, contradicting the assumption that 
$G$ is $(5,A)$-connected. 

Take a plane representation of $G$ such that the vertices in $A$ are incident with the outer face of $G$. 
Let $\sigma : V(G-a)\cup F(G-a) \longrightarrow \mathbb{Z}^+$ 
such that $\sigma(t)=4-d_{G-a}(t)$ for all $t\in V(G-a)\cup F(G-a)$. Then by Lemma~\ref{discharge}, the total charge is $\sigma(G-a)=8$.

Note that for any $t\in V(G-a)\cup F(G-a)$, if $\sigma(t) > 0$ then $t\in A$, or 
$t\in F(G-a)$ and $d_{G-a}(t) = 3$ (in which case,  $\sigma(t) = 1$). 
For each $f\in F(G-a)$ with $d_{G-a}(f)=3$, choose two vertices of $G-A$ incident with $f$ (which exists as $A$ is independent), 
and send a charge of $1/2$ from $f$ to each of these two vertices. Let $\tau$ denote the resulting charge function. 

Then $\tau(f)\leq 0$ for all $f\in F(G-a)$. Suppose $K_4^- \not\subseteq G-a$. Then, for each $v\in V(G-a)$,  
$\tau(v)$ has the upper bound in Lemma~\ref{discharge}. 
So  $\tau(v)\leq 0$ when $v\notin N(a)\cup A$, $\tau(v)\le 1$ for $v\in N(a)$ (as $A$  is independent), and $\tau(v)=\sigma(v)\leq 3$ for $v\in A-\{a\}$.

Let $f_{\infty}$ denote the outer face of $G$.  
Since $A$ is independent in $G$,  $\tau(f_{\infty})=\sigma(f_{\infty}) \le 4-(|N(a)|+7)=-3-|N(a)|$. 
If at least two vertices in $A-\{a\}$ each have degree 2 or more in $G-a$, then   $\sum\limits_{v\in A-\{a\}}\tau(v)\le 2+2+3+3=10$; so 
$$\tau(f_{\infty})+\sum\limits_{v\in A-\{a\}}\tau(v)\le -3-|N(a)|+10=7-|N(a)|.$$
If at most one vertex in   $A-\{a\}$ has degree 2 or more in $G-a$ then, since $G$ is $(5,A)$-connected, 
$d_{G-a}(f_{\infty})\ge |N(a)|+9$ and $\sum\limits_{v\in A-\{a\}}\tau(v)\le 3+3+3+3=12$; so 
$$\tau(f_{\infty})+\sum\limits_{v\in A-\{a\}}\tau(v)\le 4-(|N(a)|+9)+12=7-|N(a)|.$$

Therefore, 
$$\tau(G-a)\le \tau(f_{\infty})+\sum\limits_{v\in A-\{a\}}\tau(v) +\sum\limits_{v\in N(a)}\tau(v) \leq (7-|N(a)|)+|N(a)|= 7<8=\sigma(G-a),$$
a contradiction. \qed

\medskip

The next result will not be used in this paper, but will be used in subsequent papers in the series;
we include it here as its proof also uses discharging. 

\begin{prop}\label{6-cut2}
Let $G$ be \textsl{a} graph, $A\subseteq V(G)$, and $a\in A$  such that $|A|=6$, $|V(G)|\geq 8$, 
$(G-a, A-\{a\})$ is planar, and $G$ is $(5, A)$-connected.
Then one of the following holds.  
\begin{itemize}
 \item[$(i)$] $G-a$ contains $K_4^-$, or $G$ contains \textsl{a} $K_4^-$ in which the degree of $a$ is $2$.
 \item[$(ii)$] $G$ has \textsl{a} $5$-separation $(G_1,G_2)$ such that $a\in V(G_1\cap G_2)$, $A\subseteq V(G_1)$, $|V(G_2)|\ge 7$, and 
$(G_2-a, V(G_1\cap G_2)-\{a\})$ is planar. 
\end{itemize}
\end{prop}

\pf We may assume that 
\begin{itemize}
\item [(1)] $G$ has no 5-separation or 6-separation $(G',G'')$ such that $a\in V(G'\cap G'')$, $A\subseteq V(G')$,  $|V(G'')|\ge 8$, 
and $(G''-a, V(G'\cap G'')-\{a\})$ is planar. 
\end{itemize}
For, suppose such a separation $(G',G'')$ does exist in $G$, Then $G''$ is $(5,V(G'\cap G''))$-connected. If $|V(G'\cap G'')|=5$ then $(ii)$ holds. 
If  $|V(G'\cap G'')|=6$ we may work with $G''$ and $V(G'\cap G'')$ instead of $G$ and $A$.  

\medskip

By (1), $A$ is independent in $G$. We may further assume that 
\begin{itemize}
\item [(2)] $|V(G-a)|\ge 8$ and each vertex in $A-\{a\}$ has at least two neighbors in $G-A$. 
\end{itemize}
First, suppose $|V(G-a)|=7$.  Let $b_1, b_2$ denote the two vertices in $V(G)-A$. 
Since $G$ is $(5, A)$-connected, $d(b_i)\geq 5$ for $i\in [2]$. Therefore, since  
$(G-a, A-\{a\})$ is planar, $b_1b_2\in E(G)$ and $|N(b_1)\cap N(b_2)\cap A|\ge 2$. 
Thus, $(i)$ holds. 

So assume $|V(G-a)|\ge 8$. Then by (1), each vertex in $A-\{a\}$ has at least two neighbors in $G-A$, completing the proof of (2). 

\medskip

We may also assume that 
\begin{itemize}
\item [(3)] $G-A$ is connected.
\end{itemize}
For, otherwise, let $C_1,C_2$ be two components of $G-A$. Since $G$ is $(5, A)$-connected, 
$|N(C_i)\cap (A-\{a\})|\ge 4$ for $i\in [2]$. Hence,  $|N(C_1)\cap N(C_2)\cap (A-\{a\})|\ge 3$, contradicting the assumption that  
$(G-a, A-\{a\})$ is planar.  This proves (3). 

\medskip 

We now apply a discharging argument to $G-a$ by taking a plan representation of $G-a$ in which the vertices in $A-\{a\}$ are incident with its outer face. 
Let $\sigma : V(G-a)\cup F(G-a) \longrightarrow \mathbb{Z}^+$ 
such that $\sigma(t)=4-d_{G-a}(t)$ for all $t\in V(G-a)\cup F(G-a)$. By (2) and (3), $G-a$ is connected. So by Lemma~\ref{discharge}, the total charge is 
$\sigma(G-a)=8$.

Note that for any $t\in V(G-a)\cup F(G-a)$, if $\sigma(t) > 0$ then $t\in A-\{a\}$, or 
$t\in F(G-a)$ and $d_{G-a}(t) = 3$ (in which case,  $\sigma(t) = 1$). 
For each $f\in F(G-a)$ with $d_{G-a}(f)=3$, we may assume that $f$ is incident with two vertices in $V(G-a)-N(a)$; for, otherwise, 
$(i)$ holds. So for each $f\in F(G-a)$ with $d_{G-a}(f)=3$, we choose two vertices from $V(G-a)-N(a)$ incident with $f$, 
and send a charge of $1/2$ from $f$ to each of these two vertices. Let $\tau$ denote the resulting charge function. 
Then $\tau(f)\leq 0$ for all $f\in F(G-a)$. We may assume that  $K_4^-\not\subseteq G-a$, as otherwise $(ii)$ holds. 
Then for each $v\in V(G-a)$, $\tau(v)$ has the upper bound in Lemma~\ref{discharge}.
Thus, $\tau(v)\leq 0$ if $v\notin A$, and $\tau(v)\le 5/2$ if $v\in A-\{a\}$ (by (2)).

Denote by $f_{\infty}$ the outer face of $G-a$.  
Since $A$ is independent in $G$, $\tau(f_{\infty})=\sigma(f_{\infty}) \le -6$. 
Therefore, $\tau(f_{\infty})+\sum\limits_{v\in A-\{a\}}\tau(v) \leq -6+ 25/2 = 13/2$. 
Hence, the total charge is $$\tau(G-a)\le \tau(f_{\infty})+\sum\limits_{v\in A-\{a\}}\tau(v)  \leq 13/2<8=\sigma(G-a).$$ 
This is a contradiction. \qed

\section{Apex separations}

In this section, we prove Theorem~\ref{apexside1}. For convenience, we introduce the following terminology.  
Let $G$ be a graph and $A\subseteq V(G)$.  We say that $(G,A)$ is 
{\it plane} if $G$ is drawn in the plane with no edge crossings such that the vertices in $A$ are incident with the outer  
face. Moreover, for $a_1,\ldots, a_k\in V(G)$, we say $(G, a_1,\ldots,a_k)$ is {\it  plane} (respectively, {\it planar}) if 
$G$ is drawn (respectively, has a drawing) in a closed disc in the plane with no edge crossings  such that  
 $a_1,\ldots, a_k$ occur on the boundary of the disc in this cyclic order (clockwise or counterclockwise). 
  
We also need  a few known results.  
The first result  is a consequence of a more general result of Seymour \cite{Se80} (with 
 equivalent versions proved  in \cite{CR79, Th80, Sh80}).

\begin{lem}\label{2path} 
Let $G$ be \textsl{a} graph and let $s_1,s_2,t_1,t_2\in V(G)$ be distinct such that $G$ is $(4,\{s_1,s_2,$ $t_1,t_2\})$-connected. Then
either $G$ contains disjoint paths from $s_1$ to $t_1$ and from
  $s_2$ to $t_2$, or $(G,s_1,s_2,t_1,t_2)$ is planar. 
\end{lem} 

The next lemma is Theorem 4.3 in \cite{MY10}, where it is used to prove that if a 5-connected nonplanar graph has a 5-separation with 
one side planar and nontrivial then it contains $TK_5$.

\begin{lem}\label{w4}
Let $G$ be \textsl{a} graph and let  $a_1,a_2, a_3,a_4,a_5\in V(G)$ be distinct such that 
$(G,a_1,a_2,a_3,$ $a_4,a_5)$ is plane. Let $A=\{a_1,a_2, a_3,a_4,a_5\}$ and suppose $G$ is
$(5,A)$-connected and $|V(G)|\ge 7$. Then
there exist $w\in V(G)- A$, \textsl{a} cycle $C_w$ 
in $(G-A)-w$, and four paths
$P_1,P_2,P_3, P_4$ from $w$ to $A$ such that
\begin{itemize}
\item [$(i)$] $V(P_i\cap P_j)=\{w\}$
for $1\le i<j\le 4$, and $|V(P_i\cap C_w)|=1$ for $i\in [4]$, and 
\item [$(ii)$] there exist $1\le i\ne j\le 4$ such that $a_1$ is an end of
  $P_i$ and $a_5$ is an end of $P_j$. 
\end{itemize}
\end{lem}

We need two results from \cite{KMY15}, which may be viewed as apex versions of Lemma~\ref{w4}. 
They are used in \cite{KMY15} to deal with 5-separations with an apex side. Lemma~\ref{apex1} is Corollary 2.11 in \cite{KMY15}, and 
Lemma~\ref{apex2} is Corollary 2.12 in \cite{KMY15}. 
 
\begin{lem}\label{apex1}
Let $G$ be \textsl{a} connected graph with $|V(G)|\ge 7$, $A\subseteq V(G)$ with $|A|=5$, and $a\in A$, such that  
$G$ is $(5,A)$-connected, $(G-a,A-\{a\})$ is plane, and
$G$ has no 5-separation $(G_1,G_2)$ with $A\subseteq G_1$
and $|V(G_2)|\ge 7$. Let $w\in N(a)$ such that $w$ is not incident with  the outer face of $G-a$.
Then
\begin{itemize}
\item [$(i)$] the vertices of $G-a$ cofacial with $w$ induce \textsl{a} cycle $C_w$ in $G-a$, and
\item [$(ii)$] $G-a$ contains paths $P_1,P_2, P_3$ from $w$ to $A-\{a\}$  such that $V(P_i\cap P_j)=\{w\}$
for $1\le i<j\le 3$, and $|V(P_i\cap C_w)|=|V(P_i)\cap A|=1$ for $i\in [3]$.
\end{itemize}
\end{lem}

\begin{lem}\label{apex2}
Let $G$ be \textsl{a} connected graph with $|V(G)|\ge 7$, $A\subseteq V(G)$ with $|A|=5$, and $a\in A$, such that  
$K_4^-\not\subseteq G-a$, $G$ is $(5,A)$-connected, $(G-a,(A-\{a\})\cup N(a))$
is plane, and $G$ has no 5-separation $(G_1,G_2)$ with $A\subseteq V(G_1)$ and $|V(G_2)|\ge 7$. Then
$G-a$ is 2-connected. Moreover, either $G$ is the graph obtained from the edge-disjoint union
of the $8$-cycle $a_1b_1a_2b_2a_3b_3a_4b_4a_1$ and the $4$-cycle $b_1b_2b_3b_4b_1$ by adding $a$ and the edges $ab_i$ for $i\in [4]$,
with  $A=\{a,a_1,a_2,a_3,a_4\}$, or
there exists $w\in V(G)-A$  such that
\begin{itemize}
\item [$(i)$] the vertices of $G-a$ cofacial with $w$ induce \textsl{a} cycle $C_w$ in $G-a$ such that  
$C_w\cap D=\emptyset$, where $D$ denotes the outer cycle of $G-a$,
\item [$(ii)$] there exist paths $P_1,P_2,P_3, P_4$ in $G$ from $w$ to $A$  such that $V(P_i\cap P_j)=\{w\}$
for $1\le i<j\le 4$, and $|V(P_i\cap C_w)|=|V(P_i)\cap A|=1$  for $i\in [4]$, and
\item [$(iii)$] either $a\notin \bigcup_{i=1}^4  V(P_i)$, or  $a\in \bigcup_{i=1}^4  V(P_i)$ and we may write
$A-\{a\}=\{a_1,a_2,a_3,a_4\}$ such that
$a\in V(P_1)$, $a_i\in V(P_i)$ for $2\le i\le 4$, and  $a_1,a_2,a_3,V(P_1\cap D),a_4$ occur on $D$ in cyclic order.
\end{itemize}
\end{lem}

\medskip

{\sl Proof of Theorem~\ref{apexside1}}. We choose such separation $(G_1,G_2)$ that $G_2$ is minimal. 
Then we may assume that $G_2$ has no 5-separation $(G_2',G_2'')$ such that $|V(G_2'\cap G_2'')|\le 5$, $V(G_1\cap G_2)\subseteq V(G_2')$, and $|V(G_2'')|\ge 7$. 
For, suppose such $(G_2',G_2'')$ does exist. Then $a\notin V(G_2'\cap G_2'')$; otherwise, $(G_1\cup G_2', G_2'')$ contradicts the choice 
of $(G_1,G_2)$ (that $G_2$ is minimal). Hence, $(G_2'', V(G_2'\cap G_2''))$ is planar; so $(ii)$ holds by Lemma~\ref{planarside}.

Let $A:=V(G_1\cap G_2)=\{a, a_1, a_2, a_3, a_4\}$ such that $(G_2-a,a_1,a_2,a_3,a_4)$ is plane.   
Let $D$ denote the outer walk of $G_2-a$;  so $a_1,a_2,a_3,a_4$ occur on $D$ in clockwise order.  
We may  assume that neither $(G_1, A)$ nor $(G_2, A)$ is planar; else $(ii)$ holds by Lemma~\ref{planarside}.

Suppose there exists some $w\in N(a)\cap V(G_2-D)$.  
Then by Lemma~\ref{apex1}, the vertices of 
$G_2-a$ cofacial with $w$ induce a cycle $C_w$ in $G_2-a$, and $G_2-a$ contains paths $P_1, P_2, P_3$ 
from $w$ to $A-\{a\}$ such that $V(P_i\cap P_j)=\{w\}$ for $1\leq i < j \leq 3$, and $|V(P_i\cap C_w)|=|V(P_i)\cap A|=1$ for $i\in [3]$. 
Without loss of generality, we may assume that $V(P_i)\cap A=\{a_i\}$ for $i\in [3]$. 
Let $y\in V(G_1-A)$. Since $G$ is 5-connected, there exist independent paths $Q_1, Q_2, Q_3, Q_4$ in $G_1$ from $y$ to $a_1, a_2, a_3, a$, respectively. 
Then $C_w\cup P_1\cup P_2\cup P_3\cup Q_1\cup Q_2\cup Q_3\cup (Q_4\cup wa)$
is a $TK_5$ in $G$ in which $a$ is not a branch vertex. So $(i)$ holds.

Thus, we may assume that $N(a)\cap V(G_2)\subseteq V(D)$. By the minimality of $G_2$, $A$ is independent in $G_2$. 
Hence $(G_2-a, (A-a)\cup (N(a)\cap V(G_2)))$ is planar.

Moreover, we may assume $|N(a)\cap V(G_2-A)|\ge 2$. For otherwise, let $a'$ be the unique neighbor
of $a$ in $G_2-A$.  Then $(G_1+\{a,aa'\},G_2-a)$ is a 5-separation in $G$ and $(G_2-a,\{a',a_1,a_2,a_3,a_4\})$ is planar. 
If  $|V(G_2-a)|\ge 7$ then $(ii)$ holds by Lemma~\ref{planarside}.  So we may assume  
$|V(G_2-a)|=6$ and let $a''$ be the only vertex in $V(G_2-a)-\{a',a_1,a_2,a_3,a_4\})$. Since $G$ is 5-connected, 
$a''$ is adjacent to all of $\{a',a_1,a_2,a_3,a_4\}$. 
But this forces $a'$ to have degree at most 4 in $G$ (as $(G_2-a, \{a',a_1,a_2,a_3,a_4\})$ is planar), a contradiction.

Suppose $(ii)$ and $(iii)$ of Theorem~\ref{apexside1} both fail.
Then by Lemma~\ref{apex2},  $G_2-a$ is $2$-connected (so $D$ is a cycle) and there exists $w\in V(G_2)-A$ such that 
\begin{itemize}
\item [$(1)$]  the vertices of $G_2-a$ cofacial with $w$ induce a cycle $C_w$ in $G_2-a$ such that $C_w\cap D=\emptyset$,
\item [$(2)$] there exist paths $P_1, P_2, P_3, P_4$ in $G_2$ from $w$ to $A$ such that $V(P_i\cap P_j)=\{w\}$ 
for $1\leq i < j\leq 4$, and $|V(P_i\cap C_w)|=|V(P_i)\cap A|=1$ for $i\in [4]$, and 
\item [$(3)$] either $a\notin \bigcup\limits_{i=1}^4V(P_i)$, or $a\in \bigcup\limits_{i=1}^4V(P_i)$ and 
we may relabel $a_1, a_2, a_3, a_4$ (if necessary) such that $a\in V(P_1)$, $a_i\in V(P_i)$ for $2\leq i\leq 4$, and $a_1, a_2, a_3, V(P_1\cap D), a_4$ occur on $D$ in cyclic order.
\end{itemize}
For convenience, let $L=C_w\cup P_1\cup P_2\cup P_3\cup P_4$ and assume, without loss of generality, that $a_1,a_2,a_3,a_4$ occur on $D$ in clockwise order.  
By the planarity of $G_2-a$, we may assume that $P_i\cap D$ is 
a path for $i\in [4]$.
\medskip

{\it Case} 1. $a\notin \bigcup\limits_{i=1}^4V(P_i)$. 

Without loss of generality, we may assume that $a_i\in V(P_i)$ for $i\in [4]$. 
If $G_1-a$ has disjoint paths $S_1, S_2$ from $a_1, a_2$ to $a_3, a_4$, respectively, then 
$L\cup S_1\cup S_2$ is a $TK_5$ in $G$ in which $a$ is not a branch vertex; so $(i)$ holds.
Hence, we may assume that such $S_1, S_2$ do not exist. Then, since $G_1$ is $(5, A)$-connected, it follows from  
Lemma~\ref{2path} that $(G_1-a, a_4,a_3,a_2,a_1)$ is plane and, hence, $G_1-A$ is connected.

\medskip

{\it Subcase} 1.1. $G_1-A$ is not $2$-connected.

If $|V(G_1-A)|=2$ then let $C_1,C_2$ denote the two disjoint 1-vertex subgraphs of $G_1-A$, and  otherwise let $(C_1,C_2)$ denote a 
1-separation in $G_1-A$  with $V(C_i-C_{3-i})\ne \emptyset$ for $i\in [2]$. Since $G$ is 5-connected and $(G_1-a, a_4,a_3,a_2,a_1)$ is plane, 
$N(a)\cap V(C_i-C_{3-i})\ne \emptyset$ and $C_i-C_{3-i}$ is connected,  
for $i\in [2]$. Without loss of generality, we may assume that 
$\{a, a_1, a_2, a_3\}\subseteq N(C_1)$ and $\{a, a_1, a_3, a_4\}\subseteq N(C_2)$. 

Suppose $N(a)\cap V(G_2-A)\subseteq V(P_1\cup P_3)$. If $N(a)\cap V(G_2-A)\subseteq V(P_i)$ for some $i\in\{1,3\}$ then 
$(G_2,A)$ is planar; so $(ii)$ holds by Lemma~\ref{planarside}. Thus, we may assume that there exist $a'\in N(a)\cap V(P_1)$ and  $a''\in N(a)\cap V(P_3)$.
Let $L'$ be obtained from $L$ by replacing $P_1,P_3$ with $wP_1a',wP_3a''$, respectively, and let $R$ be a path in $G_1-\{a,a_1,a_3\}$ 
from $a_2$ to $a_4$. Then $L'\cup a'aa''\cup  R$ is a $TK_5$ in $G$ 
in which $a$ is not a branch vertex.  So $(i)$ holds. 

Hence, we may assume that  there exists $a'\in N(a)\cap V(G_2-A)$ such that $a'\notin V(P_1\cup P_3)$. Recall that $N(a)\cap V(G_2-A) \subseteq V(D)$.
Then $D-V(P_1\cup P_3)$ has a path $Q$ which is either from 
$a'$ to $a_2$ (and disjoint from $P_4$) or from $a'$ to $a_4$ (and disjoint from $P_2$). 

First, assume that  $Q$ is from $a'$ to $a_2$. Let $P_2'$ be the path in $P_2\cup Q$ between $w$ and $a'$, and $L'$ be obtained from $L$ by replacing $P_2$ with $P_2'$.
Let $R_1$ be a path in $G_1[(C_1-C_2)+\{a_1,a_3\}]$ from $a_1$ to $a_3$, and $R_2$ be a path in $G_1[(C_2-C_1)+\{a,a_4\}]$ from 
$a$ to $a_4$.  Now $L'\cup R_1\cup (R_2\cup aa')$ 
is a $TK_5$ in $G$ in which $a$ is not a branch vertex, and $(i)$ holds. 

So we may assume that  $Q$ is from $a'$ to $a_4$. Then  let  $P_4'$ be the path in $P_4\cup Q$ between $w$ and $a'$, and $L'$ be obtained from $L$ by 
replacing $P_4$ with $P_4'$. Let $R_1$ be a path in $G_1[(C_1-C_2)+\{a,a_2\}]$ from $a$ to $a_2$, and 
$R_2$ be a path in $G_1[(C_2-C_1)+\{a_1,a_3\}]$ from $a_1$ to $a_3$. Now $L'\cup (R_1\cup aa')\cup R_2$ 
is a $TK_5$ in $G$ in which $a$ is not a branch vertex, and $(i)$ holds.

\medskip
{\it Subcase} 1.2.  $G_1-A$ is $2$-connected. 

Let $D'$ denote the outer cycle of $G_1-A$, and let 
$a_i', a_i''\in N(a_i)\cap V(D')$ such that $a_1'$, $a_1''$, $a_2'$,  $a_2''$, $a_3'$, $a_3''$, $a_4'$, $a_4''$ 
occur on $D'$ in counterclockwise order, with $a_i''D'a_i'$  maximal for $i\in [4]$.
Recall that $|N(a)\cap V(D)|\ge 2$. Let $a',a''\in N(a)\cap V(D)$ be distinct.

Suppose $a',a''$ may be chosen such that  $D$ contains disjoint paths $U, W$ from $a', a''$ to $a_1, a_3$ (or $a_2,a_4$), respectively, 
and  disjoint from $P_2\cup P_4$ (or $P_1\cup P_3$). By symmetry, we may assume that $U,W$ are paths in $D$ 
 from $a', a''$ to $a_1, a_3$, respectively, and  disjoint from $P_2\cup P_4$.  
Let $L'$ be obtained  from $L$ by replacing $P_1$ with the path in $P_1\cup U$ from $w$ to $a'$ and replacing $P_3$ with the path in $P_3\cup W$ 
from $w$ to $a''$. Let $R$  be a path in $G_1-\{a,a_1,a_3\}$ from $a_2$ to $a_4$. 
Now $L'\cup a'aa''\cup R$ is a $TK_5$ in $G$ in which $a$ is not a branch vertex, and $(i)$ holds.  

So we may assume that no such $a',a''$ can be chosen. Then,  without loss of generality, we may assume that 
$D$ has disjoint paths $U, W$ from $a', a''$ to $a_3, a_4$, respectively,  
and disjoint from $P_{1}\cup P_{2}$.

Suppose $a$ has a neighbor on $a_2''D'a_4'-\{a_4', a_2''\}$ or $a_3''D'a_1'-\{a_1', a_3''\}$. By symmetry, we may assume 
the former.  Then $G_1-a_3$ has disjoint paths $R_1, R_2$ from $a,a_2$ to $a_1,a_4$, 
respectively. We modify $L$ to obtain $L'$ by replacing $P_3$ with the path in $P_3\cup U$ from $w$ to $a'$. Then  
$L'\cup (a'a\cup R_1)\cup R_2$ is a $TK_5$ in $G$ in which $a$ is not a branch vertex, and $(i)$ holds.

Thus, we may assume that $N(a)\cap V(D')\subseteq V(a_4'D'a_3'')$.  Since $(G_1, A)$ is not planar,   
$a$ has a neighbor in $V(G_1-A)-V(D')$, say $b$. 
Since $G$ is $5$-connected, $G_1$ has no $2$-cut $\{c_1, c_2\}\in V(a_4'D'a_2'')$ separating $b$ from $a_1$. 
So $(G_1-\{a,a_2,a_3,a_4\}) - a_4'D'a_2''$ has a path $R$ from $b$ to $a_1$. 
Let $L'$  be obtained from $L$ by replacing $P_3$ with the path in $P_3\cup U$ from $w$ to $a'$. Then $L'\cup (R\cup baa')\cup (a_4a_4'\cup a_4'D'a_2''\cup a_2''a_2)$ 
is a $TK_5$ in $G$ in which $a$ is not a branch vertex, and $(i)$ holds.

\medskip

{\it Case} 2. $a\in \bigcup\limits_{i=1}^4 V(P_i)$.

Recall the  notation in $(3)$. 
If $G_1-a_1$ contains disjoint paths $S_1, S_2$ from $a, a_3$ to $a_2, a_4$, respectively, 
then $L\cup S_1\cup S_2$ is a $TK_5$ in $G$ in which $a$ is not a branch vertex (so $(i)$ holds). 
Hence, we may assume that such  $S_1, S_2$ do not exist. Then, since $G$ is 5-connected, it follows from 
Lemma~\ref{2path} that $(G_1-a_1, a,a_3,a_2,a_4)$ is plane and, hence, $G_1-A$ is connected. 
We may assume that $a,a_3,a_2,a_4$ occur on the outer walk of $G_1-a_1$ in clockwise order. 

Recall that $P_i\cap D$ is a path for $i\in [4]$. So let $P_1\cap D=a'P_1a''$  and $P_i\cap D=a_i'P_ia_i$ for $2\leq i \leq 4$. 
Moreover, we may assume that $a'$ is the only neighbor of $a$ on $a'P_1a''$ (otherwise we can shorten $P_1$). 
Let $b, c\in V(D)$ such that $a_2 \in V(bDc)$ and $bDc\cap (P_3\cup P_4)=\emptyset$  and, subject to this, $bDc$ is maximal. 

Suppose there exists $v\in N(a)\cap V(bDc)$.  
Let $P_2'$ be the path in $P_2\cup bDc$ from $w$ to $v$, and let $L'$ be obtained from $L$ by replacing $P_2$ with $P_2'$. 
Let  $Q$ be a path  in $G_1-\{a,a_1,a_2\}$ from $a_3$ to $a_4$.   Then $L'\cup vaa'\cup Q$ is a $TK_5$ in $G$ in which $a$ is not a branch vertex, 
and $(i)$ holds.

So we may assume $N(a)\cap V(bDc)=\emptyset$.
However, since $(G_2, A)$ is not planar, $a$ must have a neighbor in $V(D)-V(a_3Da_4)$. Hence, by the maximality of $bDc$, one of the following holds:
\begin{itemize}
 \item [$(a)$] $a_3, a_4, a_4'$ occur on $D$ in clockwise order and there exists $v\in N(a)\cap V(a_4Da_4'-a_4)$; or
 \item [$(b)$] $a_3', a_3, a_4$ occur on $D$ in clockwise order and there exists $v\in N(a)\cap V(a_3'Da_3-a_3)$.
\end{itemize}
We consider two cases according to whether or not  $G_1-A$ is 2-connected. 

\medskip

{\it Subcase} 2.1. $G_1-A$ is $2$-connected. 

Let $D'$ denote the outer cycle of $G_1-A$ and $g_2, g_3, g, g_4$  
be neighbors of $a_2, a_3, a, a_4$ on $D'$, respectively. We may choose $g_2, g_3, g, g_4$ so that they are pairwise distinct and occur on $D'$ in 
counterclockwise order.

If $(a)$ holds, let $R_1=a_2g_2\cup g_2D'g_4\cup g_4a_4$, $R_2=ag\cup gD'g_3\cup g_3a_3$, and 
$L_a$ be obtained from $L$ by replacing $P_1$ with $wP_1a''\cup a''Da_4$  and replacing $P_4$ with $wP_4a_4'\cup vDa_4'\cup va$. 
Then $L_a\cup R_1\cup R_2$ is a   $TK_5$ in $G$ in which $a$ is not a branch vertex, and $(i)$ holds.

If $(b)$ holds, then let $R_1=a_3g_3\cup g_3D'g_2\cup g_2a_2$, $R_2=a_4g_4\cup g_4D'g\cup ga$, and
$L_b$ be obtained from $L$ by replacing $P_1$ with $wP_1a''\cup a_3Da''$ and replacing $P_3$ with $wP_3a_3'\cup a_3'Dv\cup va$. 
Then $L_b\cup R_1\cup R_2$ is a   $TK_5$ in $G$ in which $a$ is not a branch vertex, and $(i)$ holds.

\medskip

{\it Subcase} 2.2. $G_1-A$ is not $2$-connected. 

If $|V(G_1-A)|=2$ then let $C_1,C_2$ denote the two disjoint 1-vertex subgraphs of $G_1-A$. If $|V(G_1-A)|\ge 3$ then let $(C_1,C_2)$ 
be a 1-separation of $G_1-A$.  Since $G_1$ is $(5, A)$-connected and $(G_1-a_1, a, a_3,a_2,a_4)$ is planar, 
we may assume without loss of generality that 
$\{a_1, a, a_3, a_4\}\subseteq N(C_1)$ and $\{a_1, a_2, a_3, a_4\}\subseteq N(C_2)$, or $\{a_1, a_2, a_3, a\}\subseteq N(C_1)$ and 
$\{a_1, a_2, a_4, a\}\subseteq N(C_2)$.  Moreover, $C_i-C_{3-i}$ is connected for $i\in [2]$.

Suppose $\{a_1, a, a_3, a_4\}\subseteq N(C_1)$ and $\{a_1, a_2, a_3, a_4\}\subseteq N(C_2)$. Let $R_1$ be a path in 
$G[(C_1-C_2)+\{a,a_1\}]\cup a_1Da_2\cup P_2$ from $w$ to $a$, and   
$L'$ be obtained from $L$ by replacing $P_2$ with $R_1$. Let  $R_2$ be 
a path in $G_1[(C_2-C_1)+\{a_3,a_4\}]$ from $a_3$ to $a_4$. 
Then $L'\cup R_2$ is a $TK_5$ in $G$ in which $a$ is not a branch vertex, and $(i)$ holds. 
 
Thus, we may assume that $\{a_1, a_2, a_3, a\}\subseteq N(C_1)$ and $\{a_1, a_2, a_4, a\}\subseteq N(C_2)$. 

If $(a)$ holds, let $R_1$ be a path in $G_1[(C_1-C_2)+\{a,a_3\}]$ from $a$ to $a_3$, $R_2$ be 
a path in $G_1[(C_2-C_1)+\{a_2,a_4\}]$ from $a_2$ to $a_4$,  and $L_a$ be obtained from $L$ by replacing 
$P_1$ with $wP_1a''\cup a''Da_4$ 
and replacing $P_4$ with $wP_4a_4'\cup vDa_4'\cup va$.
Then $L_a\cup R_1\cup R_2$ is a   $TK_5$ in $G$ in which $a$ is not a branch vertex, and $(i)$ holds. 

So assume that $(b)$ holds. Let $R_1$ be a path in $G[(C_1-C_2)+\{a_3,a_2\}]$ from $a_3$ to $a_2$, $R_2$ be 
a path in $G[(C_2-C_1)+\{a,a_4\}]$ from $a$ to $a_4$,  and $L_b$ be obtained from $L$ by replacing $P_1$ with $wP_1a''\cup a_3Da''$ and 
replacing $P_3$ with $wP_3a_3'\cup a_3'Dv\cup va$. 
Then $L_b\cup R_1\cup R_2$ is a   $TK_5$ in $G$ in which $a$ is not a branch vertex, and $(i)$ holds. \qed

\section{Apex graphs} 

In this section, we consider apex graphs. The purpose is to give an indication how case (c) in Section 1 
will be taken care of later (by combining with later results 
on $K_4^-$). First, we give a consequence of Lemma~\ref{planarside}. 

\begin{prop}\label{contraction}
Let $G$ be \textsl{a} $5$-connected nonplanar graph, $a\in V(G)$, $T\subseteq G$ such that $a\in V(T)$, $T\cong K_2$ or $T\cong K_3$,   
$G/T$ is $5$-connected and  planar. Then $G-V(T)$ contains $K_4^-$.
\end{prop}

\pf Let $z$ denote the vertex of $G/T$ representing the contraction of $T$. Choose $A\subseteq V(G/T)$ such that $|A|=5$ and $z\in A$. 
Since $G/T$ is 5-connected and planar, $|V(G/T)|\ge 7$. Hence, by Lemma~\ref{planarside}, $G/T-z$ contains $K_4^-$. So  $G-V(T)$ contains $K_4^-$. \qed

\medskip

Next we give a result which strengthens  Corollary 2.9  in \cite{KMY15}.

\begin{prop} \label{apexvertex}
Let $G$ be \textsl{a} $5$-connected nonplanar graph and $a\in V(G)$ such that $G-a$ is planar. Then one of the following holds. 
\begin{itemize}
 \item[$(i)$] $G$ contains \textsl{a} $TK_5$ in which $a$ is not \textsl{a} branch vertex.
 \item[$(ii)$]  $G-a$ contains $K_4^-$.
 \item[$(iii)$] $G$ has \textsl{a} $5$-separation $(G_1, G_2)$ such that $V(G_1\cap G_2)=\{a, a_1, a_2, a_3, a_4\}$ and $G_2$ is the graph 
obtained from the edge-disjoint union of the $8$-cycle $a_1b_1a_2\allowbreak b_2a_3b_3a_4b_4a_1$ and the $4$-cycle $b_1b_2b_3b_4b_1$ by adding $a$ and the edges $ab_i$ for
$i\in [4]$.
\end{itemize}
\end{prop}

\pf We may assume that $G-a$ is embedded in the plane with no edge crossings. We may further assume that $K_4^-\not\subseteq G-a$; as otherwise $(ii)$ holds. 
Let $S:=\{v\in V(G-a): d_{G-a}(v)=4\}$. 

We now show $|S|\ge 8$ by applying to $G-a$ the same discharging argument used in Section 2. 
Let $\sigma : V(G-a)\cup F(G-a) \longrightarrow \mathbb{Z}^+$ 
such that $\sigma(t)=4-d_{G-a}(t)$ for all $t\in V(G-a)\cup F(G-a)$. 
Then by Lemma~\ref{discharge}, the total charge is $\sigma(G-a)=8$. 
Note that for any $t\in V(G-a)\cup F(G-a)$, if $\sigma(t) > 0$ then $t\in F(G-a)$, $d_{G-a}(t) = 3$, and $\sigma(t) = 1$
(i.e., $t$ is a triangle in $G-a$).
For each $f\in F(G-a)$ with $d_{G-a}(f)=3$, pick two of its incident vertices, and send a charge of $1/2$ from $f$ 
to each of these two vertices. Let $\tau$ denote the resulting charge function. Then $\tau(f)\leq 0$ for all $f\in F(G-a)$. 
Since we assume $K_4^-\not\subseteq G$,  for $v\in V(G-a)$, $\tau(v)$ satisfies the upper bound in Lemma~\ref{discharge}. 
Hence, $\tau(v)\le 0$ for $v\in V(G-a)-S$, and $\tau(v)\le 1$ for $v\in S$. Therefore, $|S|\geq 8$ as $\tau(G-a)=\sigma(G-a)$.

Let $b\in S$ be arbitrary, let $N(b)=\{a, b_1, b_2, b_3, b_4\}$, and let $D_b$ denote the facial cycle of $(G-a)-b$ 
containing $\{b_1,b_2,b_3,b_4\}$. Note that $G-b$ is $(5,\{a,b_1,b_2,b_3,b_4\})$-connected and $((G-b)-a, \{b_1,b_2,b_3, b_4\})$ is planar. 

Suppose  there exists $w\in N(a)-V(D_b)$.
Let $G' = (G -b) - E(G[\{a, b_1,b_2,b_3,b_4\}])$. 
We wish to apply Lemma~\ref{apex1} to $G'$ and $\{a,b_1,b_2,b_3,b_4\}$. Since $|S|\ge 8$,  $|V(G')|\geq 7$. 
Suppose $G'$ has a $5$-separation $(G'_1, G'_2)$ with $\{a, b_1,b_2,b_3,b_4\}\subseteq V(G'_1)$ and $|V(G'_2)|\geq 7$. 
Then $G$ has a separation $(G_1,G_2)$ such that $V(G_1\cap G_2)=V(G_1'\cap G_2')$, $V(G_1)=V(G_1')\cup \{b\}$ and $G_2=G_2'$.
Note that $|V(G_1')|\ge 6$ as $\{a, b_1,b_2,b_3,b_4\}$ is independent in $G'$. So $|V(G_1)|\ge 7$. 
If $a\in V(G_1\cap G_2)$ then 
we may apply Theorem~\ref{apexside1} to $G$  to conclude that $(i)$ or $(ii)$ or $(iii)$ holds. If $a\notin V(G_1\cap G_2)$ then $(G_2,V(G_1\cap G_2))$ is 
planar; so $(ii)$ holds by Lemma~\ref{planarside}. Hence,  we may assume that no such 5-separation $(G_1',G_2')$ exists in $G'$. 
Then by Lemma~\ref{apex1}, the vertices of $G'-a$ cofacial with $w$ induce a cycle $C_w$ (in $G'-a$) and $G'-a$ 
contains three paths $P_1, P_2, P_3$ from $w$ to $\{b_1,b_2,b_3,b_4\}$ such that $V(P_i\cap P_j)=\{w\}$ for $1\leq i < j\leq 3$, 
and $|V(P_i\cap C_w)|=|V(P_i)\cap \{b_1,b_2,b_3,b_4\}|=1$ for $i\in [3]$. Without loss of generality, 
assume that $b_i\in V(P_i)$ for $i\in [3]$. Now $C_w\cup baw\cup (P_1\cup b_1b)\cup (P_2\cup b_2b)\cup (P_3\cup b_3b)$ 
is a $TK_5$ in $G$ in which $a$ is not a branch vertex, and $(i)$ holds.

So we may assume that $N(a)-\{b\}\subseteq V(D_b)$.  Without loss of generality, we may assume that $b_1$, $b_2$, $b_3$, $b_4$ 
occur on $D_b$ in clockwise order. 

Since $|S|\geq 8$, we may further assume that there exists 
$v\in S\cap V(b_1D_bb_2-\{b_1, b_2\})$. Let $N(v) = \{a,v_1,v_2,v_3,v_4\}$ with $v_1 \in V(b_1D_bv)$ and $v_2 \in V(vD_bb_2)$. 
Let $D_v$ be the facial cycle of $(G-a)-v$ containing $\{v_1,v_2,v_3,v_4\}$. 
Since $b\in S$ is chosen arbitrary, we also have 
 $N(a) - \{v\} \subseteq V(D_v)$. Thus, $N(a) - \{b,v\} \subseteq V(D_v \cap D_b)=V(b_1D_bb_2)$ (since $G$ is $5$-connected).  
However, this implies that $G$ is planar, a contradiction.   \qed

\section{Topological $H$ and 5-separations}

In this section, we prove Theorem~\ref{5cut_triangle}. We fix the notation $H$ to represent the tree on six vertices, 
two of which are adjacent and of degree 3. We need a result on $TH$ in the proof of Theorem~\ref{5cut_triangle} (and also in the proof of 
Theorem~\ref{N(x1)capA} in Section 6).

Let $G$ be a graph and $u_1,u_2,a_1,a_2,a_3,a_4$ be distinct vertices of $G$.  We say that a $TH$ in $G$ is {\it rooted} at $u_1,u_2,\{a_1,a_2,a_3,a_4\}$ 
if in the $TH$, $u_1,u_2$ are of degree 3 and $a_1,a_2,a_3,a_4$ are of degree 1. 
For convenience, we use {\it quadruple} to denote $(G,u_1,u_2,A)$ where $u_1,u_2$ are distinct vertices of the graph 
$G$, $A\subseteq V(G)-\{u_1,u_2\}$, and $|A|=4$. We say that
$(G,u_1,u_2,A)$ is {\it feasible} if $G$ has a  $TH$ rooted at $u_1,u_2,A$.

In \cite{MXY15}, infeasible quadruples are characterized in terms of obstructions. Here we include the descriptions and figures of obstructions given in \cite{MXY15}. 
A quadruple $(G, u_1,u_2,A)$ is an {\it obstruction} if $G$ has subgraphs
$U_1, U_2$ (called {\it sides}) and $A_1, \ldots, A_k$ (called {\it middle parts}), where $2\le k\le 4$, such that
\begin{itemize}
\item $V(G)=V(U_1)\cup V(U_2)\cup A_{[k]}$, where $A_{[k]}=\cup_{i\in [k]} V(A_i)$,
\item $V(A_i), i\in [k],$ are pairwise disjoint,
\item $E(G-A)$ is the disjoint union of $E(U_1-A), E(U_2-A)$, and $E(A_i-A)$ for $i\in [k]$,
\item $V(U_1\cap U_2)\subseteq A\subseteq A_{[k]}$, $u_1\in V(U_1)-A_{[k]}$, and $u_2\in V(U_2)-A_{[k]}$,
\item for any $i\in [k], V(A_i)\cap A\neq \emptyset$ and $|V(A_i)\cap A|\le |V(A_i)\cap V(U_1\cup U_2)|\le |V(A_i)\cap A|+1$,
\item if $|V(A_i)|\ge 2$, then $V(A_i)\cap V(U_1\cup U_2)\cap A=\emptyset$ and $N(V(A_i)\cap A)\subseteq V(A_i)$.
\end{itemize} 
Note that $V(A_i)\cap V(U_1\cup U_2)\cap A\neq \emptyset$ iff 
$|V(A_i)|=1$, in which case $V(A_i)\subseteq A\cap V(U_1\cap U_2)$ and there is no restriction on $N(A_i)$.

\begin{figure}[!htp]
\begin{center}
\psfrag{A1}{$A_1$}
\psfrag{A2}{$A_2$}
\psfrag{A3}{$A_3$}
\psfrag{U1}{$U_1$}
\psfrag{U2}{$U_2$}
\psfrag{u1}{$u_1$}
\psfrag{u2}{$u_2$}
\psfrag{I}{type I}
\psfrag{II}{type II}
\includegraphics[scale=0.6]{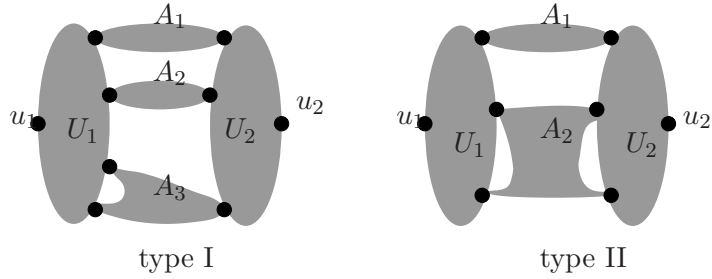}
\caption{\label{typeI-II} Obstructions of type I and type II}
\end{center}
\end{figure}

An obstruction $(G, u_1,u_2,A)$ is said to be  of  {\it type I} if  
$k=3$, $|V(A_i)\cap A|=1$ for $i\in [2]$, $|V(A_3)\cap A|=2$, 
$|V(U_i\cap A_j)|=1$ for $(i,j)\ne (1,3)$, and  $|V(U_1\cap A_3)|=2$. 

An obstruction $(G, u_1,u_2,A)$ is said to be  
of {\it type II} if $k=2$, $|V(A_1)\cap A|=1$, $|V(A_2)\cap A|=3$,
$|V(U_i\cap A_1)|=1$  for $i\in [2]$, and $|V(U_i\cap A_2)|=2$  for $i\in [2]$. 

\begin{figure}[!htp]
\begin{center}
\psfrag{A1}{$A_1$}
\psfrag{A2}{$A_2$}
\psfrag{A3}{$A_3$}
\psfrag{A4}{$A_4$}
\psfrag{U1}{$U_1$}
\psfrag{U2}{$U_2$}
\psfrag{u1}{$u_1$}
\psfrag{u2}{$u_2$}
\psfrag{III}{type III}
\psfrag{IV}{type IV}
\includegraphics[scale=0.6]{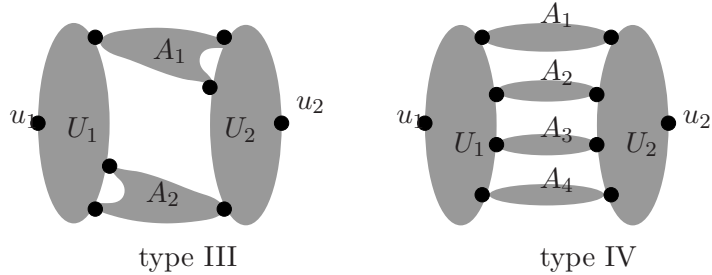}
\caption{\label{typeIII-IV} Obstructions of types III and  IV}
\end{center}
\end{figure}

An obstruction $(G, u_1,u_2,A)$ is said to be of {\it type III} if $k=2$, $|V(A_i)\cap A|=2$ for $i\in [2]$, 
$|V(U_1\cap A_1)|=|V(U_2\cap A_2)|=1$, and $|V(U_1\cap A_2)|=|V(U_2\cap A_1)|=2$. 

An obstruction $(G, u_1,u_2,A)$ is said to be
of {\it type IV} if $k=4$ and, 
for $i\in [4]$ and $j\in [2]$, $|V(A_i)\cap A|=|V(U_j\cap A_i)|=1$. 

\medskip

The following result is proved by Ma, Xie and Yu in \cite{MXY15}, which gives a characterization of feasible quadruples.

\begin{lem}
\label{H} 
Let $(G,u_1,u_2,A)$ be \textsl{a} quadruple. Then one of the following holds.
\begin{itemize}
\item [$(i)$] $(G,u_1,u_2,A)$ is feasible.
\item [$(ii)$] $G$ has \textsl{a} separation $(K,L)$ such that $|V(K\cap L)|\le 2$ and for some $i\in [2]$, 
$u_i\in V(K-L)$  and $A\cup \{u_{3-i}\}\subseteq V(L)$.
\item [$(iii)$] $G$ has \textsl{a} separation $(K,L)$ such that $|V(K\cap L)|\le 4$, $u_1,u_2\in V(K-L)$, and $A\subseteq V(L)$.
\item [$(iv)$] $(G,u_1,u_2,A)$ is an obstruction of type I, or type II, or type III, or type IV.
\end{itemize}
\end{lem}

It is easy to verify that  if $G$ is $(4,A)$-connected and  
has no separation as in $(iii)$ of Lemma~\ref{H}, then $(G,u_1,u_2,A)$ is not feasible iff it is an obstruction of type IV. Moreover, in the case when 
$G$ is $(4,A)$-connected and $(G,u_1,u_2,A)$ is an obstruction of type IV, it is straightforward to verify that for any $i\in [2]$ and $a\in A$, 
$G$ contains four independent paths: one from $u_i$ to $a$, and three from $u_{3-i}$ to the three vertices in $A-\{a\}$, respectively.

\medskip

{\sl Proof of Theorem~\ref{5cut_triangle}}. 
Let  $G$ be a 5-connected graph and $(G_1,G_2)$ be a 5-separation in $G$. Let $V(G_1\cap G_2) = \{a,a_1,a_2,a_3,a_4\}$ such that $G[\{a,a_1,a_2\}]\cong K_3$. 
We wish to prove $(iv)$; so let $u_1,u_2,u_3\in N(a)-\{a_1,a_2\}$ be pairwise distinct and  $G':=G-\{av:v\not\in \{a_1,a_2,u_1,u_2,u_3\}\}$.
 Without loss of generality, we may assume that $u_1\in V(G_2)-V(G_1)$.

\medskip

{\it Claim} 1. We may assume $u_2,u_3\in V(G_2)-V(G_1)$. 

Since $G$ is 5-connected,  $G_2-a$ has independent paths $P_1,P_2,P_3,P_4$ 
from $u_1$ to $a_1,a_2,a_3,a_4$, respectively.

First, suppose $u_2 \notin V(G_2)$; the case when $u_3 \notin V(G_2)$ can be treated in the same way. 
Since $G$ is $5$-connected,  $G_1-a$ has independent paths 
$Q_1,Q_2,Q_3,Q_4$ from $u_2$ to $a_1,a_2,a_3,a_4$, respectively.  
Then $G[\{a,a_1,a_2\}] \cup u_1a \cup P_1 \cup P_2 \cup u_2a \cup Q_1 \cup Q_2 \cup (P_3\cup Q_3)$ 
is a $TK_5$ in $G'$ with branch vertices $a,a_1,a_2, u_1,u_2$.

Next, assume that  $u_2 \in \{a_3,a_4\}$; the case when  $u_3 \in \{a_3,a_4\}$ is the same. 
Without loss of generality, we may assume that $u_2 = a_3$. 
  Let $b\in V(G_1)-V(G_2)$. 
Since $G$ is 5-connected, $G_1-a$ has independent paths $Q_1,Q_2,Q_3,Q_4$ 
from $b$ to $a_1,a_2,a_3,a_4$, respectively. Now $G[\{a,a_1,a_2\}] \cup u_1a \cup P_1\cup P_2 
\cup Q_1\cup Q_2 \cup (Q_3\cup a_3a)\cup (P_4\cup Q_4)$ is  a $TK_5$ in  $G'$ 
with branch vertices $a,a_1,a_2,b,u_1$. Thus we have Claim 1. 

\medskip

Before proceeding with the proof, we introduce additional notation. 
If there is a 4-separation $(K,L)$ in $G_2-a$ such that $u_1,u_2\in V(K-L)$ and $A:=\{a_1,a_2,a_3,a_4\}\subseteq V(L)$, we choose 
such a separation that $K$ is minimal. If such a 4-separation does not exist, we let $K=G_2$, $V(L)=A$ and $E(L)=\emptyset$. 
Let $V(K\cap L)=\{y_1,y_2,y_3,y_4\}$. Since $G$ is 5-connected, we may assume that $L$ contains disjoint paths $Y_1,Y_2,Y_3,Y_4$ from $a_1,a_2,a_3,a_4$ 
to $y_1,y_2,y_3,y_4$, respectively (by relabeling if necessary). 
\medskip

{\it Claim} 2. We may assume that $(K,u_1,u_2, \{y_1,y_2,y_3,y_4\})$ is not feasible. 

For, suppose $K$ has a $TH$ rooted at $u_1,u_2,\{y_1,y_2,y_3,y_4\}$, say $T$. 
First, assume that  there exists $i\in [2]$ such that $T-u_{3-i}$ has independent paths from $u_i$ to $y_1,y_2$, respectively. 
If $G_1-a$ has two disjoint paths $P,Q$  from $\{a_3,a_4\}$ to $\{a_1,a_2\}$ 
then $G[\{a,a_1,a_2,u_1,u_2\}] \cup T \cup Y_1\cup Y_2\cup Y_3\cup Y_4\cup P\cup Q$ is a $TK_5$ in   $G'$
with branch vertices $a,a_1,a_2,u_1,u_2$. So assume such paths $P,Q$ do not exist in $G_1-a$. Then $G_1$ has a separation $(G_1',G_2')$ 
such that $a\in V(G_1'\cap G_2')$, $|V(G_1'\cap G_2')|\le 2$, $\{a_1,a_2\}\subseteq V(G_1')$, and $\{a_3,a_4\}\subseteq V(G_2')$. 
Since $G$ is 5-connected, $|V(G_i')|=4$ for $i\in [2]$. This forces $|V(G_1)|=6$, a contradiction.

Thus,  we may assume without loss of generality that $T-u_2$ has independent paths from $u_1$ to $y_1,y_4$, respectively. 
If $G_1-a$ has disjoint paths $P_1, P_2$ from $a_3,a_4$ to $a_1,a_2$, respectively, 
then $G[\{a,a_1,a_2,u_1,u_2\}] \cup T \cup Y_1\cup Y_2\cup Y_3\cup Y_4\cup P_1 \cup P_2$ is  
a $TK_5$ in  $G'$, with branch vertices $a,a_1,a_2,u_1,u_2$. 
Hence, we may  assume that  such $P_1,P_2$ do not exist. Then by Lemma~\ref{2path},  $(G_1-a, \{a_1,a_2,a_3,a_4\})$ is planar. 
Since $|V(G_i)|\ge 7$ for $i\in [2]$, it follows from Theorem~\ref{apexside1} that $(i)$, or $(ii)$, or $(iii)$ holds, completing the proof of Claim 2.

\medskip
Since $G$ is 5-connected, it follows from  the choice of $(K,L)$, Claim 2 and Lemma \ref{H} that $(K,u_1,u_2,\{y_1,y_2,y_3,y_4\})$ is an obstruction 
of type IV. Thus, $K$ has four independent paths $R_1,R_2,R_3,R_4$, with $R_1,R_2,R_3$ from $u_1$ to  $y_1,y_2,y_3$, respectively, and 
$R_4$ from $u_2$ to $y_4$. Let $b \in V(G_1)-V(G_2)$. Since $G$ is 5-connected, $G_1-a$ contains independent paths $Q_1,Q_2,Q_3,Q_4$ from 
$b$ to $a_1,a_2,a_3,a_4$, respectively.  
Therefore, $G[\{a,a_1,a_2, u_1\}] \cup (R_1\cup Y_1)\cup (R_2\cup Y_2)\cup (R_3\cup Y_3\cup Q_3)\cup Q_1\cup Q_2\cup  (au_2\cup R_4\cup Y_4\cup Q_4)$ 
is a $TK_5$ in $G'$  with branch vertices $a,a_1,a_2, b,u_1$. Hence, $(iv)$ holds. 
\qed

\section{$6$-Separations}

In this section, we prove Theorem~\ref{N(x1)capA}.  We need the following result of Perfect \cite{Pe68}.

\begin{lem} \label{Perfect} 
Let $G$ be \textsl{a} graph, $u\in V(G)$, and $A\subseteq V(G-u)$. Suppose there exist $k$ independent paths from $u$ to distinct $a_1,\ldots, a_k\in A$, 
respectively, and otherwise disjoint from $A$. Then for any $n\ge k$, if there exist $n$
independent paths  $P_1,\ldots, P_n$ in $G$ from $u$ to $n$ distinct vertices in $A$ and otherwise disjoint from $A$ then  
$P_1,\ldots, P_n$ may be chosen so that $a_i\in V(P_i)$ for $i\in [k]$. 
\end{lem}
 
We also need the next lemma, which is due to Watkins and Mesner \cite{WM67}.
\begin{lem}
\label{Watkins}
 Let $R$ be \textsl{a} $2$-connected graph and let $y_1, y_2, y_3\in V(R)$ be distinct. Then there is no cycle through $y_1, y_2, y_3$ in $R$ 
if, and only if, one of the following statements holds. 
\begin{itemize}
\item [$(i)$] There exists \textsl{a} 2-cut $S$ in $R$ and there exist pairwise
  disjoint subgraphs $D_{y_i}\subseteq R - S$, $i\in [3]$, such
  that $y_i\in V(D_{y_i})$ and  each $D_{y_i}$ is \textsl{a} union of components of $R - S$. 
\item [$(ii)$]  There exist 2-cuts $S_{y_i}$ of $R$ and
  pairwise disjoint  subgraphs $D_{y_i}$ of $R$, $i\in [3]$, such that $y_i\in
  V(D_{y_i})$, each $D_{y_i}$ is \textsl{a} union of components of $R-S_{y_i}$,
  $S_{y_1} \cap S_{y_2} \cap S_{y_3} = \{z\}$, and $S_{y_1} - \{z\}, S_{y_2} - \{z\}, S_{y_3} - \{z\}$ are pairwise disjoint. 
\item [$(iii)$]  There exist pairwise disjoint $2$-cuts $S_{y_i}$ in $R$
  and pairwise disjoint subgraphs $D_{y_i}$ of $R$ such that
  $y_i \in V(D_{y_i})$, $D_{y_i}$ 
is \textsl{a} union of components of $R - S_{y_i}$, and $R - V(D_{y_1} \cup D_{y_2} \cup D_{y_3})$ 
has precisely two components, each containing exactly one vertex from $S_{y_i}$.
\end{itemize}
\end{lem}

\medskip

{\sl Proof of Theorem~\ref{N(x1)capA}}. 
Suppose $N(x)\cap V(G_1-G_2)= \emptyset$. We proceed to show that one of $(i)$--$(iv)$ holds. We may assume that
\begin{itemize}
\item [(1)] $|V(G_i)|\ge 8$ for $i\in[2]$; 
\end{itemize}
for, otherwise, $G$ contains $K_4^-$ (as $G$ is 5-connected) and, hence, $(ii)$ holds.  
Since $G$ is 5-connected, 
\begin{itemize}
\item [(2)] $N(x_i)\cap V(G_1-G_2)\ne \emptyset$ for $i\in [2]$, and  $N(\{x_1,x_2\})\cap V(G_2-G_1)\ne \emptyset$. 
\end{itemize}
Moreover, by the minimality of $G_1$, we have 
\begin{itemize}
\item [(3)] $|N(v_i)\cap V(G_1-G_2)|\ge 2$ for $i\in [3]$. 
\end{itemize}
By Theorem~\ref{5cut_triangle}, we may assume  
\begin{itemize}
\item [(4)] $N(v_i)\cap V(B)\ne \emptyset$ for $i\in [3]$ and for any component $B$ of $G_2-G_1$. 
\end{itemize}
For, suppose $B$ is a component of $G_2-G_1$ such that $N(v_i)\cap V(B)=\emptyset$ for some $i\in [3]$.  
If $|V(B)|=1$ then $G$ contains $K_4^-$ (as $G$ is 5-connected). So assume $|V(B)|\ge 2$. Let 
$B'$ be obtained from $B$ by adding $\{x,x_1,x_2,v_1,v_2,v_3\}-\{v_i\}$ and all edges of $G$ from $B$ to  
$\{x,x_1,x_2,v_1,v_2,v_3\}-\{v_i\}$. Now  we may apply Theorem~\ref{5cut_triangle} to the 
separation $(G-B, B')$. 
If $(i)$ or $(ii)$ of Theorem~\ref{5cut_triangle} holds then we 
have $(i)$ or $(ii)$. If $(iii)$ of Theorem~\ref{5cut_triangle} holds then $G$ contains $K_4^-$; so 
$(ii)$ holds. Now assume $(iv)$ of Theorem~\ref{5cut_triangle} holds. Then clearly, $(iii)$ holds. So we may assume (4).  

\medskip

Let $a \in V(G_1-G_2)$. Since $G$ is 5-connected and $N(x)\cap V(G_1-G_2)=\emptyset$, $G_1-x$ has  independent paths $P_1,P_2,P_3,P_4,P_5$ 
from $a$ to $x_1,x_2,v_1,v_2,v_3$, respectively.

If $N(x) = \{x_1,x_2,v_1,v_2,v_3\}$, then let $b \in V(G_2-G_1)$. Since $G$ is 5-connected, $G_2-x$ has  independent paths $Q_1,Q_2,Q_3,Q_4,Q_5$ 
from $b$ to $x_1,x_2,v_1,v_2,v_3$, respectively. Now $G[\{x,x_1,x_2\}] \cup P_1 \cup P_2 \cup (P_3\cup v_1x) \cup Q_1 
\cup Q_2 \cup (Q_4\cup v_2x) \cup (P_5\cup Q_5)$ is $TK_5$ in $G'$ with branch vertices $a,b, x,x_1,x_2$. 
Thus $(ii)$ holds with $x_3=v_1$, and $\{y_1,y_2\}=\{v_2,v_3\}$ (as
$d(x)= 5$ in this case).
So we may assume that 
\begin{itemize}
\item [(5)] $|N(x)\cap V(G_2-G_1)|\ge 1$, and let $x_3\in N(x)\cap V(G_2-G_1)$. 
\end{itemize}
We wish to prove that $(iii)$ holds; so let $y_1,y_2\in N(x)-\{x_1,x_2,x_3\}$ and $G':=G-\{xv:v\notin \{x_1,x_2,x_3,y_1,y_2\}\}$.

\medskip

{\it Case} 1.   $N(x_i)\cap V(G_2-G_1)\ne \emptyset$ for $i\in [2]$. 

\medskip 
{\it Claim} 1. We may assume $\{y_1,y_2\}\not\subseteq \{v_1,v_2,v_3\}$.

For, suppose $\{y_1,y_2\}\subseteq \{v_1,v_2,v_3\}$. Without loss of generality, let $y_i=v_i$ for $i\in [2]$.  
Let $X$ denote the component of $G_2-G_1$ containing $x_3$.  

First, suppose $N(x_i)\cap V(X)\ne \emptyset$ for $i\in [2]$. 
Then $G[X+\{x_1,x_2\}]$ contains a path $P$ from $x_1$ to $x_2$, and a path from $x_3$ to $b\in V(P)-\{x_1,x_2\}$ and internally disjoint from $P$;
so $G[X+\{x_1,x_2\}]$ has independent paths from $b$ to $x_1,x_2,x_3$, respectively. 
By Lemma~\ref{Perfect},  $G[X+\{x_1,x_2,v_1,v_2,v_3\}]$ has  independent paths $Q_1,Q_2,Q_3,Q_4$ 
from $b$ to $x_1,x_2,x_3, v_j$, respectively, and internally disjoint from $\{v_1,v_2,v_3\}$. By symmetry between $v_1=y_1$ and $v_2=y_2$, 
we may assume $j\in \{2,3\}$. Let $Q=Q_4\cup P_4$ if $j=2$, and $Q=Q_4\cup P_5$ if $j=3$.
Now $G[\{x,x_1,x_2\}] \cup P_1 \cup P_2 \cup (P_3\cup v_1x) \cup Q \cup Q_1 \cup Q_2 \cup (Q_3\cup x_3x) $ is a $TK_5$ in $G'$ 
with branch vertices $a,b,x,x_1,x_2$. 

Hence, we may assume 
by symmetry that $N(x_1)\cap V(X)=\emptyset$. Since $G$ is 5-connected, $G[X+\{x_2,v_1,v_2,v_3\}]$ 
has independent  paths $Q_1,Q_2,Q_3,Q_4$ from $x_3$ to $x_2,v_1,v_2,v_3$, respectively. 
Since $N(x_1)\cap V(G_2-G_1)\ne \emptyset$, $G_2-G_1$ has a component $X'$ such that $N(x_1)\cap V(X')\ne \emptyset$. 
Then by (4), $G[X'+\{x_1,v_3\}]$ contains a path $Q$ from $x_1$ to $v_3$. So 
$G[\{x,x_1,x_2\}] \cup P_1 \cup P_2 \cup (P_3\cup v_1x) \cup (P_4\cup Q_3)\cup Q_1 \cup (Q_4\cup Q)\cup x_3x$ 
is a $TK_5$ in $G'$ with branch vertices $a,x,x_1,x_2,x_3$.

\medskip

By Claim 1, we may let $y_1\in V(G_2-G_1)$. For convenience, let $u_1=x_3$ and $u_2=y_1$.

\medskip

{\it Claim} 2. We may assume that for any $i\in [3]$,  there do not exist  $w_1,w_2\in V(G_2-G_1)$ and two disjoint paths 
$W_1,W_2$ in $G_2-G_1$ from $w_1,w_2$ to $u_1,u_2$, respectively, 
such that $(G_2-x)-((W_1-w_1)\cup (W_2-w_2))$ has a $TH$ rooted at $w_1,w_2, \{x_1,x_2\}\cup (\{v_1,v_2,v_3\} - \{v_i\})$.

Suppose otherwise and, without loss of generality, assume that there exist $w_1,w_2\in V(G_2-G_1)$, paths $W_1,W_2$ in $G_2-G_1$, and a $TH$
in $(G_2-x)-((W_1-w_1)\cup (W_2-w_2))$ rooted at 
$w_1,w_2, \{x_1,x_2,v_1,v_2\}$. Let $K$ denote the union of $W_1,W_2$ and that $TH$. 
 
We may assume that $(G_1-x)-v_3$ contains disjoint paths $X_1,X_2$ from $x_1,x_2$ to $v_1,v_2$, respectively, and disjoint paths $X_1',X_2'$ from $x_1,x_2$ to 
$v_2,v_1$, respectively. For, otherwise,
$((G_1-x)-v_3, \{x_1,x_2,v_1,v_2\})$ is planar by Lemma~\ref{2path} and, hence by (1), $(i)$ or $(ii)$ follows from Theorem~\ref{apexside1} (as 
$(iii)$ of Theorem~\ref{apexside1} implies that $G$ contains $K_4^-$). 

Now $G[\{x,x_1,x_2,u_1,u_2\}] \cup K \cup X_1\cup X_2$ or $G[\{x,x_1,x_2,u_1,u_2\}] \cup K \cup X_1'\cup X_2'$ is
a $TK_5$ in $G'$ with branch vertices $w_1,w_2,x,x_1,x_2$.  This
completes the proof of Claim 2.

\medskip

We now set up some notation. Let $G_2'=(G_2-x)+\{v_3v_2,v_3v_1\}$.
If there exists a separation $(K,L)$ in $G_2'$ such that $|V (K \cap L)| \le 4$, $u_1, u_2 \in V(K - L)$, 
and $\{v_1,v_2,x_1,x_2\} \subseteq V(L)$, we choose such $(K,L)$ that $v_3\in V(L)$ whenever possible and, subject to this,  $K$ is minimal. 
If such a separation does not exist we
 let $K=G_2'$, $V(L)=\{v_1,v_2,x_1,x_2\}$ and $E(L)=\emptyset$.   
Note that if $L\not\subseteq K$ then  $\{v_1,v_2,x_1,x_2\}\ne V(K\cap L)$; for, otherwise, $K=G_2'$ by (4), 
a contradiction. 
Since $v_3v_1,v_3v_2\in E(G_2')$, 
if $v_3 \notin V(L)$ then $v_1,v_2 \in V(K \cap L)$. 
Let $V(K \cap L) = \{s_1,s_2,s_3,s_4\}$. 

\medskip

{\it Claim} 3. If $v_3\notin V(L)$ then $L$ has disjoint paths $S_1,S_2,S_3,S_4$ from $s_1,s_2,s_3,s_4$, respectively, to $\{v_1,v_2,x_1,x_2\}$; 
and if $v_3\in V(L)$ then $L$ has  disjoint paths $S_1,S_2,S_3,S_4$ from $s_1,s_2,s_3,s_4$, respectively, to $\{v_1,v_2,x_1,x_2\}$, or 
$\{v_1,v_3,x_1,x_2\}$, or $\{v_2,v_3,x_1,x_2\}$, and internally disjoint from $\{v_1,v_2,v_3\}$.

Suppose $L$ does not contain four disjoint paths from $\{s_1,s_2,s_3,s_4\}$ to $\{v_1,v_2,x_1,x_2\}$. Then 
 $L$ has a separation $(L_1,L_2)$ such that $|V(L_1\cap L_2)|\le 3$,  $\{v_1,v_2,x_1,x_2\}\subseteq V(L_1)$ and  $\{s_1,s_2,s_3,s_4\}\subseteq V(L_2)$. 
Moreover,  $|V(L_1\cap
L_2)|=3$ and $v_3\notin V(L_1)$, as otherwise $\{x,v_3\}\cup
V(L_1\cap L_2)$ 
would be a cut in $G$ of size at most 4.   Hence, since $v_3v_1,v_3v_2\in E(G_2')$, we have 
$v_1,v_2 \in V(L_1\cap L_2)$. 

First, assume  $v_3\notin V(L_2)$. Then $v_3\in V(K-L)$; but then $V(L_1\cap L_2)\cup \{v_3\}$ is a 
cut in $G$ of size 4 and  separating $\{u_1,u_2\}$ from $\{x_1,x_2,v_1,v_2,v_3\}$, contradicting the choice of $(K,L)$ (that $v_3\in V(L)$ whenever possible).  

So $v_3\in V(L_2)$. Let $V(L_1\cap L_2) = \{v_1,v_2,t\}$. 
Since $G$ is 5-connected, $L_2$ has four disjoint paths $S_1',S_2',S_3',S_4'$ from  $s_1,s_2,s_3,s_4$, respectively, to $\{t,v_1,v_2,v_3\}$. 
If $L_1-v_2$ has disjoint paths $X_1,X_2$ from  $\{t,v_1\}$ to $\{x_1,x_2\}$, then $S_1',S_2',S_3',S_4'$ and $X_1,X_2$ form 
four disjoint paths in $L$ from $\{s_1,s_2,s_3,s_4\}$ to $\{v_2,v_3,x_1,x_2\}$.
So assume that such $X_1,X_2$ do not exist. 
Then $L_1$ has a separation $(L_{11},L_{12})$ such that $v_2\in V(L_{11}\cap L_{12})$, $|V(L_{11}\cap L_{12})|\le 2$, $\{x_1,x_2\}\subseteq 
V(L_{11})$ and $\{t,v_1\}\subseteq V(L_{12})$.  Since $N(x_i)\cap
V(G_2-G_1)\ne \emptyset$ for $i\in [2]$, $V(L_{11}\cap L_{12})\ne
\{v_2\}$. If $|V(L_{11}-v_2)|=3$ then $G[V(L_{11}-v_2)\cup
\{x\}]$ contains $K_4^-$, and $(ii)$ holds. So we may assume
$|V(L_{11}-v_2)|\ge 4$.  Then $V(L_{11}\cap L_{12})\cup \{x,x_1,x_2\}$
is a $5$-cut in $G$. Hence, $G$ has a 5-separation $(H_1,H_2)$ such that $V(H_1\cap H_2)=V(L_{11}\cap L_{12})\cup \{x,x_1,x_2\}$, 
$L_{11}\subseteq H_1$, $L_{12}\cup L_2\cup K\cup G_1\subseteq H_2$. Clearly, $|V(H_2)|\ge 7$ and $|V(H_1)|\ge 6$. If $|V(H_1)|=6$ then 
$G$ contains $K_4^-$ (as $G$ is 5-connected); so $(ii)$ holds. If $|V(H_1)|\ge 7$ then the assertion follows from
Theorem~\ref{5cut_triangle}. This completes the proof of Claim 3.

\medskip

Note that a $TH$ in $K$ rooted at $u_1,u_2,\{s_1,s_2,s_3,s_4\}$ and the paths in Claim 3 would give a contradiction to Claim 2. Hence, 
$(K,u_1,u_2,\{s_1,s_2,s_3,s_4\})$  is not feasible. So by the minimality of $K$ and the 5-connectedness of $G$, it follows from 
 Lemma~\ref{H} that 
$(K,u_1,u_2,\{s_1,s_2,s_3,s_4\})$ is an obstruction of type I, or type II, or type III, or type IV. 
 Let $U_1,U_2,A_i,A_{[k]}$ be defined as in Section 4 with
 $\{s_1,s_2,s_3,s_4\}$ as $A$.

\medskip

{\it Claim} 4. We may assume   $(K,u_1,u_2,\{s_1,s_2,s_3,s_4\})$  is an obstruction of type IV.

If $(K,u_1,u_2,\{s_1,s_2,s_3,s_4\})$ is an obstruction of type II or type III, then by symmetry we may assume that $v_3 \notin V(U_2)-A_{[k]}$; 
now $V(U_2 \cap A_{[k]})\cup \{x\}$ is a $4$-cut in $G$, a contradiction. 

So assume  $(K,u_1,u_2,\{s_1,s_2,s_3,s_4\})$ is an obstruction of type I. 
Then $v_3\in V(U_2)-A_{[k]}$, as otherwise $V(U_2\cap A_{[k]})\cup \{x\}$ would be a 4-cut in $G$. 
Since $v_3v_1,v_3v_2\in E(G_2')$, we may assume $V(A_1) = \{v_1\}$ and $V(A_2) = \{v_2\}$. Without loss of generality, 
assume $s_1=v_1$ and $s_2=v_2$. 
Let $V(A_3\cap U_2)=\{t_{2}\}$ and  $V(A_3 \cap U_1)=\{t_{11},t_{12}\}$. 

We may assume that $A_3$ contains two disjoint paths $X_1,X_2$ from $\{s_3,s_4\}$ to  $\{t_{11},t_{12},t_2\}$. 
For,  otherwise, $A_3$ has a separation $(A_{31},A_{32})$ such that $|V(A_{31}\cap A_{32})|\le 1$, 
$\{s_3,s_4\}\subseteq V(A_{31})$, and $\{t_{11},t_{12},t_2\}\subseteq V(A_{32})$. Then  
$V(A_{31}\cap A_{32})\cup \{v_1,v_2,x\}$ is a cut in $G$ of size at most 4, a contradiction. 

 Since $G$ is 5-connected, $U_1$ has independent paths $Q^1_1,Q^1_2,Q^1_3,Q^1_4$ from 
$u_1$ to $v_1,v_2,t_{11}, t_{12}$, respectively, and $U_2$ has four independent paths $Q^2_1,Q^2_2,Q^2_3,Q^2_4$ 
from $u_2$ to $v_1,v_2,v_3, t_2$, respectively. 

If both $X_1,X_2$ end in $\{t_{11}, t_{12}\}$ then $(Q^1_1\cup Q^2_1)\cup (Q^1_3\cup Q^1_4\cup X_1
\cup X_2) \cup Q^2_2\cup Q^2_3 \cup S_3 \cup S_4 $ is a $TH$ in $G_2-x$  rooted at $u_1,u_2,\{x_1,x_2,v_2,v_3\}$, contradicting Claim 2. 
So $t_{11}\in V(X_1)$ and $t_2\in V(X_2)$. Then $(Q^1_1\cup Q^2_1)\cup Q^1_2\cup (Q^1_3\cup X_1) \cup 
Q^2_3\cup (Q^2_4\cup X_2) \cup S_3 \cup S_4$ is a $TH$ in $G_2-x$ rooted at $u_1,u_2,\{x_1,x_2,v_2,v_3\}$,
contradicting Claim 2 and completing the proof of Claim 4. 

\medskip

By Claim 4, we may assume  $s_i\in
V(A_i)$ for $i\in [4]$.  Let $V(A_i\cap U_1)=\{r_i\}$ and $V(A_i\cap U_2)=\{t_i\}$ for $i\in [4]$.

\medskip

{\it Subcase} 1.1.  $v_3\in V(K-L)$.

Then $v_1,v_2\in \{s_1,s_2,s_3,s_4\}$, and without loss of generality, assume $v_1=s_3$ and $v_2=s_4$.
Since $v_3v_1,v_3v_2\in E(G_2')$, we may assume by symmetry that $v_3\in V(A_3\cup U_2)$. Then $V(A_4)=\{v_2\}$ and $r_4=s_4=t_4=v_2$.
Since  $G$ is 5-connected, $U_1$ has independent paths $R_1,R_2,R_3,R_4$ from $u_1$ to $r_1,r_2,r_3,r_4$, respectively.
Let $X_1$ be a path in $A_1$ from $r_1$ to $s_1$, 
and $X_2$ be a path in $A_2$ from $t_2$ to $s_2$.

Suppose $v_3\in V(A_3)$.   Then $V(A_3)\ne \{v_1\}$. We may assume that $A_3$ has disjoint paths $R,R'$ from $r_3,t_3$, 
respectively, to $\{v_1,v_3\}$. For otherwise, 
$A_3$ has a separation $(A_{31},A_{32})$ such that $|V(A_{31}\cap A_{32})|\le 1$, $\{r_3,t_3\}\subseteq V(A_{31})$ and $\{v_1,v_3\}\subseteq V(A_{32})$. 
Then $V(A_{31}\cap A_{32})\cup \{s_1,s_2,s_4\}$ is a cut in $K$ separating
$\{u_1,u_2\}$ from $\{s_1,s_2,s_3,s_4\}$, contradicting the minimality
of $K$.  Since  $G$ is 5-connected, 
$U_2$ has  independent paths $T_1,T_2,T_3,T_4$ from $u_2$ to $t_1,t_2,t_3,t_4$, respectively. 
Then $(R_4\cup T_4)\cup (R_1\cup X_1)\cup (R_3\cup R)\cup (T_2\cup X_2)\cup (T_3\cup R')$ 
is a $TH$ in $K$ rooted at $u_1,u_2,\{s_1,s_2,v_1,v_3\}$. This $TH$ and the paths $S_1,S_2$ in Claim 3 (from $s_1,s_2$ to $\{x_1,x_2\}$)
form a $TH$ in $G_2-x$ rooted at  $u_1,u_2,\{x_1,x_2,v_1,v_3\}$, contradicting Claim 2.

Thus, $v_3\in V(U_2-A_3)$. Then $V(A_3)=\{v_1\}$ and $r_3=t_3=v_1=s_3$. 
 Suppose $U_2-\{t_1,t_2,t_3,t_4\}$ contains a path from $u_2$ to $v_3$.
Since $G$ is 5-connected, $U_2$ has four independent paths from $u_2$
to $\{v_3,t_1,t_2,t_3,t_4\}$ with only $u_2$ in common. So by
Lemma~\ref{Perfect},  $U_2$ contains three independent paths $Q_1,Q_2,Q_3$ from $u_2$ to $t_p,t_q,v_3$, respectively, 
such that $p\in \{1,2\}$, $q\in \{3,4\}$, and 
disjoint from $\{t_1,t_2,t_3,t_4\}-\{t_p,t_q\}$.  Without loss of generality, we may assume $p=2$ and $q=3$. 
Then $(R_3\cup Q_2)\cup Q_3\cup (Q_1\cup X_2)\cup (R_1 \cup X_1)\cup R_4$ is a $TH$ in $K$ rooted at $u_1,u_2, \{s_1,s_2,v_2,v_3\}$. 
 This $TH$ and the paths $S_1,S_2$ in Claim 3 (from $s_1,s_2$ to $\{x_1,x_2\}$) form a  $TH$ in $G_2-x$ 
rooted at  $u_1,u_2,\{x_1,x_2,v_2,v_3\}$, contradicting Claim 2.

Hence, $U_2-\{t_1,t_2,t_3,t_4\}$ contains no path from $u_2$ to $v_3$. Then $U_2$ has a separation $(U_{21},U_{22})$ such that 
$V(U_{21}\cap U_{22})=\{t_1,t_2,t_3,t_4\}$, $v_3\in V(U_{21})$ and $u_2\in V(U_{22})$. By the minimality of $K$, $K-\{s_1,s_2,s_3,s_4\}$ must have a path from $v_3$ to $u_2$. 
Thus we may assume that  $t_1\ne s_1$ and $G[V(U_{21})]-\{t_2,t_3,t_4\}$ contains a path
$Q$ from $v_3$ to $t_1$. Let $Q_1, Q_2,Q_3,Q_4$ be independent paths in $U_{22}$ from $u_2$ to $t_1,t_2,t_3,t_4$, respectively. 
Then $(Q_3\cup R_3)\cup (Q_1\cup Q)\cup (Q_2\cup X_2)\cup (R_1\cup  X_1)\cup R_4$ is a $TH$ in $K$ rooted at $u_1,u_2,\{s_1,s_2,v_2,v_3\}$. This $TH$ 
and the paths $S_1,S_2$ in Claim 3 (from $s_1,s_2$ to $\{x_1,x_2\}$) give a $TH$ in $G_2-x$  
rooted at  $u_1,u_2,\{x_1,x_2,v_2,v_3\}$, contradicting Claim 2.

\medskip

{\it Subcase} 1.2.  $v_3\notin V(K-L)$. 

Since $G$ is 5-connected, $|V(A_i)|=1$ or $|V(A_i)|=3$; as otherwise $\{r_i,s_i,t_i,x\}$ would be a cut in $G$ of size 4. 
Moreover, by the choice of $(K,L)$ (the minimality of $K$), 
if $|V(A_i)|=3$ then $r_is_it_i$ is a path in $G$.  
Thus, since $G$ is 5-connected, $K$ contains eight independent paths $R_1,R_2,R_3,R_4$ from $u_1$ to $s_1,s_2,s_3,s_4$, respectively, and 
$Q_1,Q_2,Q_3,Q_4$ from $u_{2}$ to $s_1,s_2,s_3,s_4$, respectively.  Let $S_1,S_2,S_3,S_4$ be the paths in Claim 3 with $s_i\in S_i$ for $i\in [4]$. 

Suppose for some $i\in [4]$ and $j\in [2]$, $s_i$ has at least two neighbors in $U_j$. Without loss of generality, assume $i=j=1$. Then, $r_1=s_1=t_1$. 
Since $G$ is 5-connected, $U_1-\{r_2,r_3,r_4\}$ has a path from $s_1$
to $u_1$.  Moreover, $U_1$ has two independent paths from $s_1$ to $\{u_1,r_2,r_3,r_4\}$ with only $s_1$ in common. For, otherwise, 
$U_1$ has a separation $(U_{11},U_{12})$ such that $|V(U_{11}\cap U_{12})|\le 1$, $s_1\in V(U_{11}-U_{12})$, and 
$\{u_1,r_2,r_3,r_4\}\subseteq V(U_{12})$. Now $|V(U_{11})|\ge 3$, as $s_1$ has at least two neighbors in $U_1$. Hence, $V(U_{11}\cap U_{12})\cup \{x,s_1\}$
is a cut in $G$ of size at most 3, a contradiction.  So by Lemma~\ref{Perfect}, $U_1$ has two independent paths
$R_1',R_2'$ from $s_1$ to $u_1, r_p$, respectively, with $p\in \{2,3,4\}$, and internally 
disjoint from  $\{r_2, r_3,r_4\}$. Without loss of generality, we
may assume $p=2$. Now $Q_1\cup (Q_3\cup S_3)\cup (Q_4\cup S_4)\cup  S_1 \cup (R_2'\cup
r_2s_2\cup S_2)$ form a $TH$ in $(G_2-x) -(R_1'-s_1)$ rooted at
$s_1,u_2,\{x_1,x_2,v_1,v_2\}$, or $s_1,u_2,\{x_1,x_2,v_1,v_3\}$, or
$s_1,u_2,\{x_1,x_2,v_2,v_3\}$, However, this  contradicts Claim 2. 

Thus, we may assume that no such $i,j$ exist. Then each $s_i$ has at least two neighbors in $L-\{s_1,s_2,s_3,s_4\}$, unless $s_i\in \{x_1,x_2,v_1,v_2,v_3\}$.

Suppose there exists some $s_i\notin \{x_1,x_2,v_1,v_2,v_3\}$, say $i=1$. 
Let $t$ be the end of
$S_1$ other than $s_1$.  By the similar argument as above, $L$ has two
independent paths from $s_1$ to $\{t,v_1,v_2,v_3\}\cup V(S_2\cup S_3\cup S_4)$.
Then by Lemma~\ref{Perfect}, there exist two independent paths
$S_1',R$ from $s_1$ to $t,r$, respectively, with $r\in \{v_1,v_2,v_3\}\cup  V(S_2\cup
S_3\cup S_4)$ and internally disjoint from $\{v_1,v_2,v_3\}\cup  V(S_2\cup
S_3\cup S_4)$. If $r\in \{v_1,v_2,v_3\}-V(S_2\cup S_3\cup S_4)$ then $S_1'\cup R\cup Q_1\cup (Q_3\cup S_3)\cup (Q_4\cup S_4)$ is a 
$TH$ in $(G_2-x)-(R_1-s_1)$ rooted at  $s_1,u_2,\{v_1,v_2,x_1,x_2\}$, or
$s_1,u_2,\{v_1,v_3, x_1,x_2\}$, or $s_1,u_2,\{v_2,v_3, x_1,x_2\}$, contradicting  Claim 2. 
Hence, without loss of generality, we may assume $r\in V(S_2)$ and let $r'$ be the end of $S_2$ other than $s_2$. 
Now $S_1'\cup (R\cup rS_2r')))\cup Q_1\cup (Q_3\cup S_3)\cup (Q_4\cup S_4)$ is a  
$TH$ in $(G_2-x)-(R_1-s_1)$ rooted at  $s_1,u_2,\{v_1,v_2,x_1,x_2\}$, or
$s_1,u_2,\{v_1,v_3,x_1,x_2\}$, or $s_1,u_2,\{v_2,v_3, x_1,x_2\}$. This
contradicts  Claim 2. 

Hence, $\{s_1,s_2,s_3,s_4\}\subseteq  \{x_1,x_2,v_1,v_2,v_3\}$. By (4) and by symmetry between $x_1$ ad $x_2$, we may assume that 
$s_i=v_i$ for $i\in [3]$ and $s_4=x_2$. Since $N(x_1)\cap V(G_2-G_1)\ne \emptyset$, $N(x_1)\cap V(L-K)\ne\emptyset$. Thus, by (4),  $G[L-K+\{x_1,v_1\}]$
has a path $T$ from $v_1$ to $x_1$. If $(G_1-x)-\{v_1,v_2\}$ has two independent paths $T_1,T_2$ from $v_3$ to $x_1,x_2$, respectively then 
$G[\{x,x_1,x_2,u_1\}]\cup T_1\cup T_2\cup (Q_3\cup u_2x)\cup (R_1\cup T)\cup R_3\cup R_4$ is a $TK_5$ in $G'$ with branch vertices 
$u_1,v_3, x, x_1, x_2$. So assume such $T_1,T_2$ do not exist.  Then 
$G_1-x$ has a separation  $(G_1',G_1'')$ such that $|V(G_1'\cap G_1'')|\le 3$, $\{v_1,v_2\}\subseteq V(G_1'\cap G_1'')$, $
v_3\in V(G_1')$ and $\{x_1,x_2\}\subseteq V(G_1'')$. 
Since $N(v_3)\cap V(G_1-G_2)\ge 2$ by (3), $V(G_1'\cap G_1'')\cup \{v_3\}$ is a cut in $G$ of size at most 4, a contradiction.

\medskip

{\it Case} 2. $N(x_i)\cap V(G_2-G_1)=\emptyset$ for some $i\in [2]$.

Without loss of generality, we may assume that $N(x_1)\cap V(G_2-G_1)=\emptyset$. 

\medskip

{\it Claim} 1.  We may assume $\{y_1,y_2\}\not\subseteq \{v_1,v_2,v_3\}$.

For, otherwise, we may assume $y_i=v_i$ for $i\in [2]$.  Since $G$ is 5-connected, $G_2-\{x,x_1\}$ contains 
 independent paths $Q_1,Q_2,Q_3,Q_4$ from $x_3$ to $v_1,v_2,v_3, x_2$, respectively. 
If $G_1-\{x,x_2\}$ has a cycle $C$ containing $\{x_1,v_1,v_2\}$, then 
$G[\{v_1,v_2, x,x_1,x_3\}] \cup C\cup (Q_4\cup x_2x_1)\cup Q_1\cup Q_2$ is a $TK_5$ in $G'$ with branch vertices $v_1,v_2,x,x_1,x_3$.
So we may assume that such a cycle $C$ does not exist. We wish to apply Lemma~\ref{Watkins}. 

First, we show that $G_1-\{x,x_2\}$ is 2-connected. Note that  $G_1-G_2$ is connected by the minimality of
$G_1$. Suppose $G_1-\{x,x_2\}$ has a cut vertex, say $v$. 
Then $v\in V(G_1-G_2)$, and   $G_1-\{x,x_2\}$ has a separation 
$(G_{11},G_{12})$ such that $V(G_{11}\cap G_{12})=\{v\}$, $x_1\in V(G_{11})$ and $V(G_{12})\cap \{v_1,v_2,v_3\}\ne \emptyset$. 
But this implies that $V(G_{12})-\{v_1,v_2,v_3\}=\{v\}$ or
$\{v_1,v_2,v_3\} \subseteq V(G_{12}) $ (as $G$ is $5$-connected and $N(x)\cap V(G_1-G_2)=\emptyset$). If $V(G_{12})-\{v_1,v_2,v_3\}=\{v\}$, then any
vertex in  $V(G_{12})\cap \{v_1,v_2,v_3\}$ has at most one neighbor in $G_1-G_2$, contradicting (3).
So $\{v_1,v_2,v_3\} \subseteq V(G_{12})$. This implies $V(G_{11})=\{x_1,v\}$ (as $G$ is $5$-connected). But now $N(x_1) \subseteq \{x,x_2,v\} $, which is a contradiction.

So by Lemma~\ref{Watkins},  $G_1-\{x,x_2\}$ has a 2-cut $\{z_1,z_2\}$ separating $x_1$ from $\{v_1,v_2\}$. Let $D_{x_1}$ 
denote the component of 
$(G_1-\{x,x_2\})-\{z_1,z_2\}$ containing $x_1$. Then $|V(D_{x_1}) - \{x_1, v_3\}|\le 1$; 
otherwise $S:=\{x,x_1,x_2,z_1,z_2,v_3\}$ is a $6$-cut in $G$ such that $G-S$ has a component strictly contained in $G_1-G_2$, contradicting  the
choice of $(G_1,G_2)$ (that $G_1$ is minimal).
Suppose $|V(D_{x_1}) - \{x_1, v_3\}|= 1$, and let $v \in V(D_{x_1}) - \{x_1, v_3\}$. Since $G$ is 5-connected and $N(x)\cap V(G_1-G_2)=\emptyset$, 
$N(v)=\{v_3,x_1,x_2,z_1,z_2\}$. So $G[\{v,x,x_1,x_2\}]\cong K_4^-$, and $(ii)$ holds. Therefore, we may assume $V(D_{x_1}) \subseteq  \{x_1, v_3\}$. 
Since $N(x_1)\cap V(G_2-G_1)=\emptyset$, 
$N(x_1) = \{x,x_2,z_1,z_2,v_3\}$. Hence, $v_3\in V(D_{x_1})$. Since $|N(v_3)\cap V(G_1-G_2)|\ge 2$,  we have 
$z_1,z_2 \in N(v_3)$. Therefore, $G[\{x_1,v_3,z_1,z_2\}]$ contains $K_4^-$, and $(ii)$ holds. 

\medskip

By Claim 1, let $y_1\in V(G_2-G_1)$. For convenience, let $u_1:=x_3$ and $u_2:=y_1$. 

\medskip

{\it Claim} 2. We may assume that there exist $w_1,w_2\in V(G_2-G_1)$ and two disjoint paths $W_1,W_2$ in $G_2-G_1$ from $w_1,w_2$ to $u_1,u_2$, respectively, 
such that $(G_2-\{x,x_1\})-((W_1-w_1)\cup (W_2-w_2))$ has a $TH$ rooted at $w_1,w_2,\{x_2,v_1,v_2,v_3\}$.

First, we set up some notation. 
If $G_2-\{x,x_1\}$ has a separation $(K,L)$ such that $|V(K\cap L)|\le 4 $, $u_1,u_2\in V(K-L)$ and $\{x_2,v_1,v_2,v_3\}\subseteq V(L)$ we select 
$(K,L)$ so that $K$ is minimal; in this case, since $G$ is 5-connected, $|V(K\cap L)|=4$, and we let $V(K\cap L)=\{s_1,s_2,s_3,s_4\}$. 
If such a separation $(K,L)$ does not exist in $G_2-\{x,x_1\}$, we let $K=G_2-\{x,x_1\}$,  $V(L)=\{s_1,s_2,s_3,s_4\}=\{x_2,v_1,v_2,v_3\}$, and $E(L)=\emptyset$.

Since $G$ is 5-connected, $L$ has four disjoint paths $S_1,S_2,S_3,S_4$ from $s_1,s_2,s_3,s_4$, respectively, to $\{x_2,v_1,v_2,v_3\}$.  
Thus, we may assume that $K$ has no  $TH$ rooted at $u_1,u_2,\{s_1,s_2,s_3,s_4\}$; otherwise such $TH$ and $S_1,S_2,S_3,S_4$ form a $TH$ in
$G_2-\{x,x_1\}$ rooted at $u_1,u_2,\{x_2,v_1,v_2,v_3\}$ and, hence,  Claim 2 holds by letting, for $i\in [2]$, $u_i=w_i$ and $W_i$ be the trivial path 
consisting of $w_i$ only. 

Since $G$ is 5-connected and by the choice of $(K,L)$ and Lemma~\ref{H}, $(K,u_1,u_2,\{s_1,s_2,s_3,s_4\})$ is an obstruction of type IV.
 As before, we use the notation $U_1,U_2$ and $A_i$ from Section 4, with $s_i\in A_i$ for $i\in [4]$ and with $\{s_1,s_2,s_3,s_4\}$ as $A$. 
For $i\in [4]$, let $V(A_i\cap U_1)=\{r_i\}$ and $V(A_i\cap U_2)=\{t_i\}$. 
Since $G$ is 5-connected, $|V(A_i)|=1$ or $|V(A_i)|=3$, and $r_is_it_i$ is a path in $A_i$; so $K$ has eight independent paths 
$P^i_1, P^i_2,P^i_3,P^i_4$ from $u_i$ to $s_1,s_2,s_3,s_4$,
respectively. 

Suppose $\{s_1,s_2,s_3,s_4\}=\{x_2,v_1,v_2,v_3\}$. Without loss of generality, let $s_i=v_i$ for $i\in [3]$, and $s_4=x_2$.
If  $G[G_1-\{x,x_2\}]$ has three independent paths $P_1,P_2,P_3$ from $x_1$ to $v_1,v_2,v_3$, respectively, then 
$G[\{u_1,u_2, x,x_1,x_2\}]\cup P^1_4\cup P^2_4\cup (P^1_3\cup P^2_3)\cup (P_1\cup P^1_1)
\cup (P_2\cup P^2_2)$ is a $TK_5$ in $G'$ with branch vertices $u_1,u_2,x,x_1,x_2$. 
So assume that such paths do not exist. 
Then $G[G_1-x]$ has a separation $(G_1',G_1'')$ such that 
$x_2\in V(G_1'\cap G_1'')$, $|V(G_1'\cap G_1'')|\le 3$, $x_1\in V(G_1')$ and $\{v_1,v_2,v_3\}\subseteq V(G_1'')$. Since 
$N(x_1)\cap V(G_2-G_1)=\emptyset$, $V(G_1'\cap G_1'')\cup \{x_1\}$ is a cut in $G$ of size at most $4$, a contradiction.

Thus, we may assume without loss of generality that $s_1\notin \{x_2,v_1,v_2,v_3\}$. Since $G$ is 5-connected, $s_1$ has at least two 
neighbors in $L$, or at least two neighbors in $U_i$ for some $i\in [2]$. 

First, assume that $s_1$ has at least two neighbors in some $U_i$, say $U_1$ (by symmetry).
Then  $V(A_1)=\{s_1\}$ and $s_1 \in V(U_1 \cap U_2)$.  We claim that $U_1$ has two independent paths 
from $s_1$ to two distinct  vertices in $\{u_1,r_2,r_3,r_4\}$. For, otherwise, $U_1$ has a separation $(U_{11},U_{12})$ such that $|V(U_{11}\cap U_{12})|\le 1$, 
$s_1\in V(U_{11})$ and $\{u_1,r_2,r_3,r_4\}\subseteq V(U_{12})$. Since $s_1$ has at least two neighbors in $U_1$, we see that $|V(U_{11})|\ge 3$. So $V(U_{11}\cap 
U_{12})\cup \{x,s_1\}$ is cut in $G$ of size at most 3, a contradiction. Since $U_1-\{r_2,r_3,r_4\}$ has a path from $s_1$ to $u_1$, 
it follows from Lemma~\ref{Perfect} that $U_1$ contains two independent paths $Q_1,Q_2$ from $s_1$ to $u_1,r_i$, respectively, for some $i\in 
\{2,3,4\}$, and disjoint from $\{r_2,r_3,r_4\}-\{r_i\}$. Without loss of generality, we may assume that $Q_2$ ends at $r_2$. 
Now $P^2_1\cup (P^2_3\cup S_3)\cup (P^2_4\cup S_4)\cup (Q_2\cup r_2s_2\cup S_2)\cup S_1$ is a $TH$ in $(G_2-\{x,x_1\})-(Q_1-s_1)$ rooted 
at $s_1,u_2,\{x_2,v_1,v_2,v_3\}$. Hence Claim 2 holds with $w_1=s_1$, $W_1=Q_1$, $w_2=u_2$, and $W_2 = w_2$.

So we may assume that  $s_1$ has at least two neighbors in $L$. 
Let $s\in \{x_2,v_1,v_2,v_3\}$ be the end of $S_1$ other than $s_1$. 
Since $s_1$ has at least two neighbors in $L$ and $G$ is 5-connected, $L$ has two independent 
paths $S_1',S_2'$ from $s_1$ to $\{s\}\cup V(S_2\cup S_3\cup S_4)$ and internally disjoint from $\{s\}\cup S_2\cup S_3\cup S_4$. 
Thus by Lemma~\ref{Perfect} (and the existence of $S_1$), $L$ has independent paths $S_1',S_2'$ from $s_1$ to 
$s,s'$, respectively, and internally disjoint from $S_2\cup S_3\cup S_4$, with $s'\in V(S_2\cup S_3\cup S_4)$. 
Without loss of generality, 
we may further assume that $S_2'$ ends at $s'\in V(S_2)$ and let $t$ be the end of $S_2$ in $\{x_2,v_1,v_2,v_3\}$. Now 
$S_1'\cup (s_1S_2's'\cup s'S_2t\cup P^2_1\cup (P^2_3\cup S_3)\cup (P^2_4\cup S_4)$ is a $TH$ in $(G_2-\{x,x_1\})-(P_1^1-s_1)$ rooted 
at $s_1,u_2,\{x_2,v_1,v_2,v_3\}$. Hence,  Claim 2 holds with $w_1=s_1$,  $W_1=P^1_1$,  $w_2=u_2$, and $ W_2 = w_2$. 
This completes the proof of Claim 2.

\medskip

By Claim 2, let $Z$ denote the union of $W_1,W_2$ and a $TH$ in  $(G_2-\{x,x_1\})-((W_1-w_1)\cup (W_2-w_2))$  rooted at $w_1,w_2,\{x_2,v_1,v_2,v_3\}$.
Without loss of generality, we may assume that $Z-w_2$ has independent paths from $w_1$ to $u_1,x_2,v_1$, respectively, and
$Z-w_1$ has independent paths from $w_2$ to $u_2,v_2,v_3$, respectively. 

It remains to show that $(iv)$ holds for $i=j=1$ if we cannot find the desired $TK_5$ for $(iii)$. 
Suppose $G_1-x$ contains three independent paths, with one from $x_2$ to $v_3$ 
and two from $x_1$ to $v_1,v_2$,  respectively, 
or with one from $x_2$ to $v_2$ and two from $x_1$ to $v_1,v_3$, respectively. 
It is easy to verify that these paths and $Z$ form a $TK_5$ in $G'$
with branch vertices $w_1,w_2,x,x_1,x_2$. So we may assume that any three independent paths from $\{x_1,x_2\}$ to $v_1,v_2,v_3$, respectively, with two 
from $x_1$ and one from $x_2$,  
must contain a path from $x_2$ to $v_1$. Hence $(iv)$ holds. \qed

\end{document}